\documentclass[12pt]{article}
\usepackage{latexsym}
   \usepackage{amsfonts}

\newtheorem{theorem}{Theorem}[section]
\newtheorem{lemma}[theorem]{Lemma}
\newtheorem{proposition}[theorem]{Proposition}
\newtheorem{corollary}[theorem]{Corollary}
\newtheorem{definition}[theorem]{Definition}

\newtheorem{remark}[theorem]{Remark}

\def \Q {{\mathbb Q}}
\def \R {{\mathbb R}}
\def \C {{\mathbb C}}

\def \Im { {\mathrm{Im}} }
\def \Re { {\mathrm{Re}} }

\begin{document}
\title{\LARGE{\bf{Shuffle Product\\ 
for Multiple Dedekind Zeta Values\\
over Imaginary Quadratic Fields}}}

\author{
Michael Dotzel, Ivan Horozov\\
}
\date{}

\maketitle

\begin{abstract}

Multiple Dedekind zeta values were recently defined by the second author. In a separate paper, the second author constructed double shuffle relations in some cases as a response to questions asked by Richard Hain 
and Alexander Goncharov.
 
In this paper, we develop the technique for obtaining more shuffle relations and produce many examples of shuffle products over an imaginary quadratic field. We also define the notion of self shuffle of a (multiple) Dedekind zeta value and use it at many instances. We define a refinement of the multiple Dedekind zeta values. Our key examples are self shuffles of the Dedekind zeta at 2 and at 3, the shuffle product of the Dedekind zeta of 2 times itself, and the shuffle product of the Dedekind zeta at 2 times the Dedekind zeta at 3.

We obtain one unexpected result that the self shuffle of multiple Dedekind zeta at (1,2) minus the self shuffle  of the twisted (with a permutation) multiple Dedekind zeta value at (1,2) is a very simple expression in terms of the refined multiple Dedekind zeta values.
\end{abstract}

\tableofcontents
\setcounter{section}{-1}

\section{Introduction}
Recently, the last author has defined multiple Dedekind zeta values (MDZVs) as a number field analogue of Euler's multiple zeta values (MZVs) in the same way as the Dedekind zeta values are a number field analogue of the Riemann zeta values. 

Before recalling the definitions for each type of zeta value, we would like to say a few words regarding the purpose of this paper. It is well known that a product of two MZV can be expressed as a sum of such. One way to obtain this sum is to write each MZV as an infinite series. Then one can express the product of such series as a sum of similar types of series, which are again MZVs. This process is called a {\it{stuffle product}}.
There is another way of expressing a product of MZVs as a sum of such, which produces a different formula. It is based on an integral representation of MZV as Chen's iterated integrals. For such integrals, there is a shuffle product, which produces a sum of iterated integrals. Each of the new iterated integrals is again a MZV. The second way of expressing a product of two MZVs as a sum of such using iterated integrals is called a
{\it{shuffle product}}. The shuffle and the stuffle product are often called double shuffle relations.

Many people have asked whether the MDZVs possess a double shuffle relation similar to the double shuffle for MZVs. The first answer to this question is given in \cite{double shuffle}, where only one example of double shuffles is considered in details. 

his paper is the first one of a sequence of three papers where we develop more systematic approach and produce much more examples of shuffles and stuffles of MDZVs.

The first paper of this sequence deals with shuffle products. We refine the definition so that we can perform stuffle products. Theoretically, the new ingredient is exactly the refinement of MDZVs. Simply said, the MDZVs defined in \cite{MDZF} can be expressed as a finite sum of the MDZVs that are defined in this paper. The advantage of the new definition is twofold. First, the new MDZVs have an integral representation as well as shuffle and stuffle products. The previous definition of MDZV has only an integral representation and a shuffle product, but not a stuffle product. Second, in the paper \cite{MDZF}, we prove the analytic continuation of the (old) MDZVs by performing exactly the same refinement into a finite sum, which are the new MDZVs. For each of the new MDZVs, we prove the analytic continuation separately. Therefore the analytic continuation holds for the old MDZVs.

The second paper of this sequence deals with the stuffle product. Since the examples of stuffle products and shuffle products are quite lengthy, we decided to separate them into two papers. Moreover, both the first and the second papers feature new definitions and new examples, which is another reason for separating them into two papers.

The careful reader would notice that we mentioned about a sequence of three papers regarding the double shuffle. We described briefly the first one, establishing shuffle product, and the second one, establishing the stuffle product, which are two papers not three. 

The third paper will use both the shuffle and stuffle product in order to obtain linear relations among MDZVs (over imaginary quadratic fields). It might seem at first that we will be using only the previous two papers and produce many examples. However, for MDZVs there is a new phenomenon: there exists a "self shuffle" and a "self stuffle" relation of sorts, which we will call a {\it{selfie}} and a {\it{souffle}}, respectively. 

If we perform a shuffle product of two (multiple) Dedekind zeta values, then we can express them as a sum of MDZVs. However, the double shuffle does not stop there. Each of the new MDZVs can be self stuffled (souffled). After that point neither the shuffle nor the stuffle can be applied. This process leads to one formula.

Conversely, if we take the same two (multiple) Dedekind zeta values, then we can perform a stuffle product an express them in a different way as a sum of MDZVs. Each of the new MDZVs can be self shuffled (selfied). After this point neither the shuffle nor the stuffle can be applied. This leads to another formula.

Comparing the two formulas, we obtain a linear relation among MDZVs, which is the goal of the third paper.

Let us recall the various zeta values that we mentioned. Riemann zeta values are defined by
\[\zeta(n)=\sum_{k=1}^\infty \frac{1}{k^n}.\]
Multiple zeta values (MZVs) are defined as
\[\zeta(n_1,\dots,n_d)=\sum_{0<k_1<k_2<\dots<k_d}\frac{1}{k_1^{n_1}k_2^{n_2}\dots k_d^{n_d}}.\]
For $d=1$ or $2$, they were defined by Euler \cite{Eu}.
Over a finite extension of $\Q$, there is an analogue of the Riemann zeta, called the Dedekind zeta, which is defined as follows:
\[\zeta_K(n)
=
\sum_{{\mathfrak{a}}\neq (0)}
\frac{1}{N({\mathfrak{a}})^n},
\]
where the sum is over all ideals ${\mathfrak{a}}$ different from the zero ideal of the ring of algebraic integers in the number field $K$. 

For the definition of multiple Dedekind zeta values (MDZVs) we follow the paper \cite{MDZF}.
In order to define MDZVs, first we consider the above sum, however instead of summing over all ideals different from $(0)$, we sum over all non-zero principal ideals. (If we sum over a fixed ideal class of non-principal ideals, there is trick to reduce the sum to a summation over principal ideals.) We also use a definition of discrete cones over a ring of algebraic integers. 
Given a cone $C$, we define a MDZV as
\[\zeta_{K;C}(n_1,n_2)=\sum_{\alpha,\beta\in C}\frac{1}{N(\alpha)^{n_1}N(\alpha+\beta)^{n_2}}.\]

One advantage of this definition is that it has a natural integral representation. Let $dt$ be the standard measure on the real line. Let \[DT_8=dt_1 dt_1' du_1 du_1' dt_2 dt_2' du_2 du_2'\] be the $8$-fold product of such measures. Then
\begin{eqnarray*}
&&\zeta_{K;C}(2,2)=\\
&&=\sum_{\alpha,\beta\in C}\frac{1}{N(\alpha)^{2}N(\alpha+\beta)^{2}}=\\
&&=\int_{t_1>t_1'>u_1>u_1'>0,\,\,t_2>t_2'>u_2>u_2'>0}
\exp(-\alpha_1t_1-\alpha_2t_2-\beta_1u_1-\beta_2u_2)
DT_8,
\end{eqnarray*}
where $K$ is a quadratic number field and $\alpha_1$ and $\alpha_2$ are the values of the two different complex embeddings of $K$ into $\C$; similarly for $\beta$.

The key idea in the computation of the above integral is the following lemma, which is an exercise in calculus.
\begin{lemma} If $Re(\alpha)>0$ then
\[\int_{t>u}\exp(-\alpha t)dt=\frac{\exp(-\alpha u)}{\alpha}.\]
\end{lemma}

A key part of this paper is based on iterated integrals. In the one dimensional setting, over paths, that idea was extensively examined by K.-T. Chen (see \cite{Ch1}, \cite{Ch2}). The higher dimensional analogue of iterated path integrals, we call "iterated integrals on membranes". First they were defined by the first author in \cite{ModSym}, as a first attempt to generalize Manin's construction \cite{M} of a noncommutative modular symbol to a noncommutative Hilbert modular symbol on a Hilbert modular surface. That preprint was improved and published as \cite{Hilbert}. Another idea from the preprint \cite{ModSym} was the construction of multiple Dedekind zeta values, which was completed in \cite{MDZF}. The idea of such iterated integrals was inspired by a talk by Ronny Brown on cubical categories based on the paper \cite{cubes}. 

Independently, Thomas Tradler, Scott O. Wilson, Mahmoud Zeinalian \cite{TWZ}, defined similar types of iterated integrals to examine analogues of connections in dimension 2, called gerbes.

For applications in number theory, especially Dedekind zetas, we heavily use the idea of a ``cone", which started with a paper of Don Zagier \cite{Z} for real quadratic fields, developed further by Shintani \cite{Sh} for totally real number fields and used by \cite{C} to construct $p$-adic interpolation of Dedekind zeta values. In the paper \cite{MDZF}, we have extended the notion of a cone in order to be suitable for any number ring (not only totally real). In the current paper, we use the idea of ``cone'' in order to denote ``positive cone'' from the paper \cite{MDZF}. This means that the integers in the imaginary  quadratic ring that that we are summing over have positive real part.

In the paper \cite{MDZF}, for the purpose of proving analytic continuation, based on distribution, \cite{GSh}, we constructed a refinement of MDZVs, essentially by cutting the cubes in the definition of the MDZVs into simplices in a suitable way. A similar idea was used also in the first paper on shuffles of of MDZV, \cite{double shuffle}. It is used in the current paper as well. 

Other analogues of MDZV over imaginary quadratic rings were examined by \cite{G}.


\section{Self shuffle of (multiple) Dedekind zeta values}
Let us consider the simplest case of convergent Dedekind zeta values, namely, $\zeta_K(2)$, where $K$ is an imaginary quadratic field. 

Let ${\cal{O}}_K$ be the ring of integers in the number field $K$. We define a domain $C$ in the ring ${\cal{O}}_K$
\[
C
=
\{\alpha\in {\cal{O}}_K \, | \, \Re(\alpha)>0 \}
\cup 
\{\alpha\in {\cal{O}}_K \, | \, \Re(\alpha)=0,\, \Im(\alpha)>0 \}.
\]

\subsection{Self shuffle of $\zeta_{K,C}(2)$}
We define Dedekind zeta at $2$ by
\[\zeta_K(2)
=
\sum_{\alpha\in C}
\frac{1}{N(\alpha)^2}.
\]
Note that this definition differs from the classical one. The essential difference here is that the sum is only over principal ideals. 
Let $\alpha_1=\alpha$ and $\alpha_2$ be two complex conjugates. Then
\[\zeta_K(2)=
\sum_{\alpha\in C}\int_{t_1>u_1>0, \,\, t_2>u_2>0}
\exp(-\alpha_1t_1-\alpha_2t_2)
dt_1du_1dt_2du_2\]
Here the domain of integration consists of a product of two simplices. 
Alternatively, the domain of integration consists of two independent chains of inequalities between the variables. The {\it{self-shuffle}} consists of all possible shuffles of the two chains of inequalities, namely,
\begin{eqnarray*}
&&\zeta_K(2)=\\
&&=\sum_{\alpha\in C}\int_{t_1>u_1>0, \,\, t_2>u_2>0}
\exp(-\alpha_1t_1-\alpha_2t_2)
dt_1du_1dt_2du_2=\\
&&=\sum_{\alpha\in C}
\left(\int_{t_1>u_1>t_2>u_2>0}+\int_{t_1>t_2>u_1>u_2>0}+\right.\\
&&+\int_{t_1>t_2>u_2>u_1>0}+\int_{t_2>t_1>u_1>u_2>0}+\\
&&\left.
+\int_{t_2>t_1>u_2>u_1>0}+ \int_{t_2>u_2>t_1>u_1>0}
\right)
\exp(-\alpha_1t_1-\alpha_2t_2)
dt_1du_1dt_2du_2=\\
&&=\sum_{\alpha\in C}
\left(
\frac{1}{\alpha_1^2(\alpha_1+\alpha_2)^2}
+
\frac{1}{\alpha_1^1(\alpha_1+\alpha_2)^3}
+\right.\\
&&\left.+
\frac{1}{\alpha_1^1(\alpha_1+\alpha_2)^3}
+
\frac{1}{\alpha_2^1(\alpha_1+\alpha_2)^3}
+\right.\\
&&\left.+
\frac{1}{\alpha_2^1(\alpha_1+\alpha_2)^3}
+
\frac{1}{\alpha_2^2(\alpha_1+\alpha_2)^2}
\right).
\end{eqnarray*}

\begin{definition}
\label{def self-sh of zeta(n)}
Given a cone $C$, we define the following multiple Dedekind zeta values $\zeta^{\rho}(n,m)$, where $\rho$ is a permutation of the two elements $\{1,2\}$ by
\[\zeta^{\rho}(n,m)=\sum_{\alpha_1\in C}\frac{1}{(\alpha_{\rho(1)})^m(\alpha_1+\alpha_2)^n},\]
where $\alpha_1$ and $\alpha_2$ are Galois conjugates.
\end{definition}
Then we obtain the following.
\begin{proposition}
\label{prop self-sh of zeta(2)}
(a) The self-shuffle of the Dedekind zeta $\zeta_{K,C}(2)$ is
\[\zeta_K(2)=\zeta^{(1)}(2,2)+2\zeta^{(1)}(1,3)+2\zeta^{(12)}(1,3)+\zeta^{(12)}(2,2).\]
(b) If the cone $C$ is Galois invariant then $\zeta^{(12)}(2,2)=\zeta^{(1)}(2,2)$ and $\zeta^{(12)}(1,3)=\zeta^{(1)}(1,3)$. In that case
\[\zeta_{K,C}(2)=2\zeta^{(1)}(2,2)+4\zeta^{(1)}(1,3).\]
\end{proposition}
\begin{remark}
In Proposition \ref{prop self-sh of zeta(2)} part (b), the formula resembles the following shuffle of MZV:
\[\zeta(2)\zeta(2)=2\zeta(2,2)+4\zeta(1,3).\] In both cases the (integral) shuffle process is exactly the same.
The difference is that the integrands are different, Thus the values are different.
\end{remark}

\subsection{Self shuffle of $\zeta_{K,C}(3)$}

For an integral representation of the Dedekind zeta value $\zeta_{K,C}(3)$, we use the same formula as in the paper \cite{MDZF}.

\begin{eqnarray}
\nonumber
\zeta_{K,C} (3) & = & \sum_{\alpha \in \mbox {\scriptsize{K}}} \frac{1}{N(\alpha)^3} \\
\label{zeta(3)}
& = & \sum \int_{t_1 > u_1 > v_1 >0; t_2 > u_2 > v_2 >0} \exp(-\alpha_1t_1 - \alpha_2t_2)DT_6
\end{eqnarray}

In order to obtain a refinement of this Dedekind zeta value, we need to write all the variables in the domain of integration as a single chain of strict inequalities. We have to list all possible ways of writing all the variables in a single chain of inequalities such that they satisfy each of the two chains of inequalities from the domain of integration in Formula (\ref{zeta(3)}). The following list of inequalities exhausts all such possible shuffles of variables between the two chains:

\begin{enumerate}
\item $t_1 > u_1 > v_1 > t_2 > u_2 > v_2>0$
\item $t_1 > u_1 > t_2 > v_1 > u_2 > v_2>0$
\item $t_1 > u_1 > t_2 > u_2 > v_1 > v_2>0$
\item $t_1 > u_1 > t_2 > u_2 > v_2 > v_1>0$
\item $t_1 > t_2 > u_1 > v_1 > u_2 > v_2>0$
\item $t_1 > t_2 > u_1 > u_2 > v_1 > v_2>0$
\item $t_1 > t_2 > u_1 > u_2 > v_2 > v_1>0$
\item $t_1 > t_2 > u_2 > u_1 > v_1 > v_2>0$
\item $t_1 > t_2 > u_2 > u_1 > v_2 > v_1>0$
\item $t_1 > t_2 > u_2 > v_2 > u_1 > v_1>0$
\item $t_2 > t_1 > u_1 > v_1 > u_2 > v_2>0$
\item $t_2 > t_1 > u_1 > u_2 > v_1 > v_2>0$
\item $t_2 > t_1 > u_1 > u_2 > v_2 > v_1>0$
\item $t_2 > t_1 > u_2 > u_1 > v_1 > v_2>0$
\item $t_2 > t_1 > u_2 > u_1 > v_2 > v_1>0$
\item $t_2 > t_1 > u_2 > v_2 > u_1 > v_1>0$
\item $t_2 > u_2 > t_1 > u_1 > v_1 > v_2>0$
\item $t_2 > u_2 > t_1 > u_1 > v_2 > v_1>0$
\item $t_2 > u_2 > t_1 > v_2 > u_1 > v_1>0$
\item $t_2 > u_2 > v_2 > t_1 > u_1 > v_1>0$
\end{enumerate}

Therefore, we obtain the following refinement, where the enumeration of the above 20 domains of of integration correspond to the order of the summation in the formula below.

\begin{eqnarray*}
\zeta_{K,C} (3) 
& = &\sum \int_{t_1 > u_1 > v_1 >0; t_2 > u_2 > v_2 >0} \exp(-\alpha_1t_1 - \alpha_2t_2)DT_6 \\
& = & \sum \frac{1}{\alpha_1^3 (\alpha_1 + \alpha_2)^3} + \sum \frac{1}{\alpha_1^2 (\alpha_1 + \alpha_2)^4} \\
& + & \sum \frac{1}{\alpha_1^2 (\alpha_1 + \alpha_2)^4} +  \sum \frac{1}{\alpha_1^2 (\alpha_1 + \alpha_2)^4} \\
& + & \sum \frac{1}{\alpha_1 (\alpha_1 + \alpha_2)^5} + \sum \frac{1}{\alpha_1 (\alpha_1 + \alpha_2)^5} \\        
& + & \sum \frac{1}{\alpha_1 (\alpha_1 + \alpha_2)^5} + \sum \frac{1}{\alpha_1 (\alpha_1 + \alpha_2)^5} \\ 
& + & \sum \frac{1}{\alpha_1 (\alpha_1 + \alpha_2)^5} + \sum \frac{1}{\alpha_1 (\alpha_1 + \alpha_2)^5} \\ 
& + & \sum \frac{1}{\alpha_2 (\alpha_1 + \alpha_2)^5} + \sum \frac{1}{\alpha_2 (\alpha_1 + \alpha_2)^5} \\ 
& + & \sum \frac{1}{\alpha_2 (\alpha_1 + \alpha_2)^5} + \sum \frac{1}{\alpha_2 (\alpha_1 + \alpha_2)^5} \\ 
& + & \sum \frac{1}{\alpha_2 (\alpha_1 + \alpha_2)^5} + \sum \frac{1}{\alpha_2 (\alpha_1 + \alpha_2)^5} \\ 
& + & \sum \frac{1}{\alpha_2^2 (\alpha_1 + \alpha_2)^4} + \sum \frac{1}{\alpha_2^2 (\alpha_1 + \alpha_2)^4} \\
& + & \sum \frac{1}{\alpha_2^2 (\alpha_1 + \alpha_2)^4} + \sum \frac{1}{\alpha_2^3 (\alpha_1 + \alpha_2)^3} \\
& = & \zeta^{(1)}(3,3) + 3\zeta^{(1)}(2,4) + 6\zeta^{(1)}(1,5) + 6\zeta^{(12)}(1,5) + 3\zeta^{(12)}(2,4) + \zeta^{(12)}(3,3).
\end{eqnarray*}
If the cone $C$ is Galois invariant then
\begin{proposition}
\label{prop self-sh of zeta(3)}
\begin{equation}
\label{refinement of zeta(3)}
\zeta_{K,C} (3) = 2\zeta^{(1)}(3,3) + 6\zeta^{(1)}(2,4) + 12\zeta^{(1)}(1,5).
\end{equation}
\end{proposition}

{\bf{Remark:}} Note that the refinement of the Dedekind zeta $\zeta_{K,C}(3)$ from Formula (\ref{refinement of zeta(3)}) resembles the shuffle product of the zeta values
\[
\zeta(3)\zeta(3) = 2\zeta(3,3) + 6\zeta(2,4) + 12(1,5).
\]

We leave as an exercise for the reader to prove that if the shuffle product of 
\[\zeta(n)\zeta(n)=\sum_i c_i\zeta(a_i,b_i)\]
then
\[\zeta_{K,C}(n)=\sum_i c_i\zeta^{(1)}(a_i,b_i).\]

\subsection{Definition of refinement of MDZVs}
Let $\alpha_1$ and $\alpha_3$ be elements of a cone $C$. Let $\alpha_2$ and $\alpha_4$ be the Galois conjugates of $\alpha_1$ and $\alpha_3$, respectively. We give the following refinement of the definition of MDZV.
\begin{definition}
\label{def refined MDZV}
\[\zeta^{(1)}(a,b,c,d)=
\sum_{\alpha_1,\alpha_3\in C}
\frac{1}{(\alpha_1)^a(\alpha_1+\alpha_2)^b(\alpha_1+\alpha_2+\alpha_3)^c(\alpha_1+\alpha_2+\alpha_3+\alpha_4)^d}.\]
If $\rho$ is a permutation of the four elements $\{1,2,3,4\}$, then we define
\[\zeta^{\rho}(a,b,c,d)=
\sum_{\alpha_1,\alpha_3\in C}
\frac{1}{(\alpha_{\rho(1)})^a(\alpha_{\rho(1)}+\alpha_{\rho(2)})^b(\alpha_{\rho(1)}+\alpha_{\rho(2)}+\alpha_{\rho(3)})^c(\alpha_1+\alpha_2+\alpha_3+\alpha_4)^d}.\]
\end{definition}

\subsection{Simplification using classes of permutation}
\label{subsec algorithm}
In this subsection we construct an algorithm that simplifies expressions consisting of many terms of MDZVs. The key idea of the simplification is to use classes of permutations instead of permutations that appear in the new definition of MDZVs. 

Assume that the summation of the terms in a MDZV is over $\alpha_1$ and $\alpha_3$ each of which belongs to the same cone $C$. Moreover, we assume that the cone $C$ is Galois invariant; that is, if we conjugate each element of the cone $C$ using the non-trivial element of the Galois group $Gal(K/\Q)$ then we obtain a new cone that coincides with the first cone $C$. Let $\alpha_2$ and $\alpha_4$ be the Galois conjugates of $\alpha_1$ and $\alpha_3$, respectively. Then the possible permutations in the MDZV of the above type are permutations $S_4$ of the four elements $\{1,2,3,4\}$. 

There are non-trivial permutations $\rho\in S_4$ such that $\zeta^{\rho}(a,b,c,d)=\zeta^{(1)}(a,b,c,d)$.
If $\rho$ is the Galois conjugation of the first cone then the above equality holds since the cone $C$ is Galois invariant. In terms of particular permutation $\rho_1=(12)$ since it permutes $\alpha_1$ with $\alpha_2$. Similarly, the conjugation of the second cone $\rho_2=(34)$ will permute $\alpha_3$ with $\alpha_4$. Also we can interchange the two cones. That would not change the value of the MDZV. In this case, $\rho_3=(13)(24)$. We can multiply $\rho_1$, $\rho_2$ and $\rho_3$ in any order and the resulting permutation would give the same MDZV as if we had the trivial permutation.

In is straightforward to check that the subgroup of $S_4$ generated by $\rho_1$, $\rho_2$ and $\rho_3$ consists of the following elements:
\[H=\{(1),(12), (34), (12)(34), (13)(24), (14)(23), (1324), (1423)\}.\]
Then if $\rho$ is a permutation of four elements in $S_4$ and $h\in H$ then $\zeta^{\rho h}(a,b,c,d)=\zeta^{\rho}(a,b,c,d)$, since the permutation $\rho h$ acts on the set $\{1,2,3,4\}$ by $h$, which does not change the value and then by $\rho$. Thus the cosets we need to consider are of the type $\rho H$. Since there are $24$ elements in $S_4$ and $8$ elements in $H$ then there are only three cosets in $S_4/H$.

Now we can state the algorithm for simplification of sums MDZVs. Let $g_1$ and $g_2$ be elements of $S_4$. First we check whether their are in the same $H$-coset. If $g_2=g_1 h$ for some $h\in H$ then they are in the same coset. In fact this is always true if $g_1$ and $g_2$ are in the same $H$-coset. In that case $g_1^{-1}g_2=h\in H.$

Therefore we have the algorithm
\[\zeta^{g_1}(a,b,c,d)=\zeta^{g_2}(a,b,c,d),\]
when $g_1^{-1}g_2\in H$.

\subsection{Self shuffle of $\zeta^{(1)(1)}_{K,C}(1,2)$}
The definition and the integral representation of $\zeta^{(1)(1)}_{K,C}(1,2)$ is taken from the paper \cite{MDZF} Section ``Examples". The shuffle computation below is obtained by shuffling the two chains of inequalities $t_1 > u_1 > v_1 >0$ and $t_2 > u_2 > v_2 >0$. Their shuffles gives us the new types of MDZVs as given in Definition \ref{def refined MDZV}

\begin{eqnarray*}
\zeta^{(1)(1)}_{K,C} (1,2) & = & \sum_{\alpha,\beta \in \mbox {\scriptsize{K}}} \frac{1}{N(\alpha) N(\alpha + \beta)^2} \\
& = & \sum \int_{t_1 > u_1 > v_1 >0; t_2 > u_2 > v_2 >0} \exp(-\alpha_1t_1 - \alpha_2t_2 - \beta_1u_1 - \beta_2u_2) \\
& = & \zeta^{(1)}(1,2;1,2)+ \zeta^{(1)}(1,1;2,2) + 2\zeta^{(1)}(1,1;1,3)\\
&& + \zeta^{(23)}(1,1;2,2) + 2\zeta^{(23)}(1,1;1,3) +  2\zeta^{(234)}(1,1;1,3)\\
&& + \zeta^{(234)}(1,1;2,2)+ \zeta^{(132)}(1,1;2,2) + 2\zeta^{(132)}(1,1;1,3) \\
&& +  2\zeta^{(1342)}(1,1;1,3) + \zeta^{(1342)}(1,1;2,2) + 2\zeta^{(13)(24)}(1,1;1,3)\\
&& + \zeta^{(13)(24)}(1,1;2,2) +  \zeta^{(13)(24)}(1,2;1,2).
\end{eqnarray*}

Based on the algorithm defined in Subsection \ref{subsec algorithm} we can simplify the expression if the cone  $C$ is Galois invariant. Then 
\begin{eqnarray*}
\zeta_{K,C} (1,2)
&=&  2\zeta^{(1)}(1,2;1,2) \\
&&+ 2\zeta^{(1)}(1,1;2,2) + 4\zeta^{(23)}(1,1;2,2)+\\
&&+ 4\zeta^{(1)}(1,1;1,3) +  8\zeta^{(23)}(1,1;1,3) 
\end{eqnarray*}

\subsection{Self shuffle of $\zeta^{(12)(1)}_{K,C}(1,2)$}
The definition and the integral representation of $\zeta^{(12)(1)}_{K,C}(1,2)$ is taken from the paper \cite{MDZF} Section ``Examples". We proceed with the computation exactly in the same way as in the previous Subsection.
\begin{eqnarray*}
\zeta_K^{(12),(1)} (1,2) & = & \sum_{\alpha,\beta \in \mbox {\scriptsize{K}}} \frac{1}{\alpha_2\beta_1 N(\alpha + \beta)^2} \\ 
&& =  \sum \int_{t_1 > u_1 > v_1 >0; t_2 > u_2 > v_2 >0} \exp(-\alpha_1u_1 - \alpha_2t_2 - \beta_1t_1 - \beta_2u_2) \\
& =  &\zeta^{(12)}(1,2;1,2) + \zeta^{(12)}(1,1;2,2) + 2\zeta^{(12)}(1,1;1,3) \\
&&+ \zeta^{(123)}(1,1;2,2) + 2\zeta^{(123)}(1,1;1,3)+2\zeta^{(1234)}(1,1;1,3)  \\
&& +  \zeta^{(1234)}(1,1;2,2) + \zeta^{(13)}(1,1;2,2) + 2\zeta^{(13)}(1,1;1,3) \\
&& +  2\zeta^{(134)}(1,1;1,3) + \zeta^{(134)}(1,1;2,2) + 2\zeta^{(13)(24)}(1,1;1,3) \\
&& +  \zeta^{(13)(24)}(1,1;2,2) + \zeta^{(13)(24)}(1,2;1,2).  
\end{eqnarray*}
If the cone is Galois invariant then
\begin{eqnarray*}
\zeta_K^{(12),(1)} (1,2) 
& =&  2\zeta^{(1)}(1,2;1,2)\\
&& + 2\zeta^{(1)}(1,1;2,2)+   3\zeta^{(23)}(1,1;2,2) +  \zeta^{(1234)}(1,1;2,2)\\
&&+ 4\zeta^{(1)}(1,1;1,3) + 6\zeta^{(23)}(1,1;1,3)+2\zeta^{(1234)}(1,1;1,3)  \\
\end{eqnarray*}

As a consequence we obtain the following
\begin{corollary}
\begin{eqnarray*}
\zeta_K^{(1),(1)} (1,2) - \zeta_K^{(12),(1)} (1,2) 
& =&\zeta^{(23)}(1,1;2,2) - \zeta^{(1234)}(1,1;2,2)\\
&&+ 2\zeta^{(23)}(1,1;1,3)-2\zeta^{(1234)}(1,1;1,3). 
\end{eqnarray*}
\end{corollary}

\section{Shuffle product of $\zeta_{K,C}(2)\zeta_{K,C}(2)$}
\label{sec DZ2}
\subsection{Preliminary}
For the Dedekind zeta at $2$ we have
\[\zeta_{K,C}(2)=2\zeta^{(1)}(2,2)+4\zeta^{(1)}(1,3),\]
when the cone is Galois invariant. In that case,
Therefore,
\[\zeta_{K,C}(2)\zeta_{K,C}(2)=4\zeta^{(1)}(2,2)\zeta^{(1)}(2,2)+8\zeta^{(1)}(2,2)\zeta^{(1)}(1,3)+16\zeta^{(1)}(1,3)\zeta^{(1)}(1,3).\]

In order to complete the shuffle of the product of $\zeta_{K,C}(2)$ with itself, we need to compute the following three shuffles products
$\zeta^{(1)}(2,2)\zeta^{(1)}(2,2)$, $\zeta^{(1)}(2,2)\zeta^{(1)}(1,3)$ and $\zeta^{(1)}(1,3)\zeta^{(1)}(1,3)$.
This problems are solved in the following four subsections, where the first one is dealing only with the needed inequalities among the variable.

\subsection{Needed inequalities for the shuffles}
This subsection deals with inequalities needed in the following three subsections dealing with similar types of shuffle products. More precisely, here we compute the shuffle of
\[t_1>t_2>t_3>t_4>0\]
and
\[u_1>u_2>u_3>u_4>0.\]

\begin{enumerate}

\item $t_1>t_2>t_3>t_4>u_1>u_2>u_3>u_4>0$
\item $t_1>t_2>t_3>u_1>t_4>u_2>u_3>u_4>0$
\item $t_1>t_2>t_3>u_1>u_2>t_4>u_3>u_4>0$
\item $t_1>t_2>t_3>u_1>u_2>u_3>t_4>u_4>0$
\item $t_1>t_2>t_3>u_1>u_2>u_3>u_4>t_4>0$

***

\item $t_1>t_2>u_1>t_3>t_4>u_2>u_3>u_4>0$
\item $t_1>t_2>u_1>t_3>u_2>t_4>u_3>u_4>0$
\item $t_1>t_2>u_1>t_3>u_2>u_3>t_4>u_4>0$
\item $t_1>t_2>u_1>t_3>u_2>u_3>u_4>t_4>0$

***

\item $t_1>t_2>u_1>u_2>t_3>t_4>u_3>u_4>0$
\item $t_1>t_2>u_1>u_2>t_3>u_3>t_4>u_4>0$
\item $t_1>t_2>u_1>u_2>t_3>u_3>u_4>t_4>0$

***

\item $t_1>t_2>u_1>u_2>u_3>t_3>t_4>u_4>0$
\item $t_1>t_2>u_1>u_2>u_3>t_3>u_4>t_4>0$

***

\item $t_1>t_2>u_1>u_2>u_3>u_4>t_3>t_4>0$

***
***

\item $t_1>u_1>t_2>t_3>t_4>u_2>u_3>u_4>0$
\item $t_1>u_1>t_2>t_3>u_2>t_4>u_3>u_4>0$
\item $t_1>u_1>t_2>t_3>u_2>u_3>t_4>u_4>0$
\item $t_1>u_1>t_2>t_3>u_2>u_3>u_4>t_4>0$

***

\item $t_1>u_1>t_2>u_2>t_3>t_4>u_3>u_4>0$
\item $t_1>u_1>t_2>u_2>t_3>u_3>t_4>u_4>0$
\item $t_1>u_1>t_2>u_2>t_3>u_3>u_4>t_4>0$

***

\item $t_1>u_1>t_2>u_2>u_3>t_3>t_4>u_4>0$
\item $t_1>u_1>t_2>u_2>u_3>t_3>u_4>t_4>0$

***

\item $t_1>u_1>t_2>u_2>u_3>u_4>t_3>t_4>0$

***
***

\item $t_1>u_1>u_2>t_2>t_3>t_4>u_3>u_4>0$
\item $t_1>u_1>u_2>t_2>t_3>u_3>t_4>u_4>0$
\item $t_1>u_1>u_2>t_2>t_3>u_3>u_4>t_4>0$

***

\item $t_1>u_1>u_2>t_2>u_3>t_3>t_4>u_4>0$
\item $t_1>u_1>u_2>t_2>u_3>t_3>u_4>t_4>0$

***

\item $t_1>u_1>u_2>t_2>u_3>u_4>t_3>t_4>0$

***
***

\item $t_1>u_1>u_2>u_3>t_2>t_3>t_4>u_4>0$
\item $t_1>u_1>u_2>u_3>t_2>t_3>u_4>t_4>0$

***

\item $t_1>u_1>u_2>u_3>t_2>u_4>t_3>t_4>0$

***
***

\item $t_1>u_1>u_2>u_3>u_4>t_2>t_3>t_4>0$

***
***
***

\item $u_1>t_1>t_2>t_3>t_4>u_2>u_3>u_4>0$
\item $u_1>t_1>t_2>t_3>u_2>t_4>u_3>u_4>0$
\item $u_1>t_1>t_2>t_3>u_2>u_3>t_4>u_4>0$
\item $u_1>t_1>t_2>t_3>u_2>u_3>u_4>t_4>0$

***

\item $u_1>t_1>t_2>u_2>t_3>t_4>u_3>u_4>0$
\item $u_1>t_1>t_2>u_2>t_3>u_3>t_4>u_4>0$
\item $u_1>t_1>t_2>u_2>t_3>u_3>u_4>t_4>0$

***

\item $u_1>t_1>t_2>u_2>u_3>t_3>t_4>u_4>0$
\item $u_1>t_1>t_2>u_2>u_3>t_3>u_4>t_4>0$

***

\item $u_1>t_1>t_2>u_2>u_3>u_4>t_3>t_4>0$

***
***

\item $u_1>t_1>u_2>t_2>t_3>t_4>u_3>u_4>0$
\item $u_1>t_1>u_2>t_2>t_3>u_3>t_4>u_4>0$
\item $u_1>t_1>u_2>t_2>t_3>u_3>u_4>t_4>0$

***

\item $u_1>t_1>u_2>t_2>u_3>t_3>t_4>u_4>0$
\item $u_1>t_1>u_2>t_2>u_3>t_3>u_4>t_4>0$

***

\item $u_1>t_1>u_2>t_2>u_3>u_4>t_3>t_4>0$

***
***

\item $u_1>t_1>u_2>u_3>t_2>t_3>t_4>u_4>0$
\item $u_1>t_1>u_2>u_3>t_2>t_3>u_4>t_4>0$

***

\item $u_1>t_1>u_2>u_3>t_2>u_4>t_3>t_4>0$

***
***

\item $u_1>t_1>u_2>u_3>u_4>t_2>t_3>t_4>0$

***
***
***

\item $u_1>u_2>t_1>t_2>t_3>t_4>u_3>u_4>0$
\item $u_1>u_2>t_1>t_2>t_3>u_3>t_4>u_4>0$
\item $u_1>u_2>t_1>t_2>t_3>u_3>u_4>t_4>0$

***

\item $u_1>u_2>t_1>t_2>u_3>t_3>t_4>u_4>0$
\item $u_1>u_2>t_1>t_2>u_3>t_3>u_4>t_4>0$

***

\item $u_1>u_2>t_1>t_2>u_3>u_4>t_3>t_4>0$

***
***

\item $u_1>u_2>t_1>u_3>t_2>t_3>t_4>u_4>0$
\item $u_1>u_2>t_1>u_3>t_2>t_3>u_4>t_4>0$

***

\item $u_1>u_2>t_1>u_3>t_2>u_4>t_3>t_4>0$

***
***

\item $u_1>u_2>t_1>u_3>u_4>t_2>t_3>t_4>0$

***
***
***

\item $u_1>u_2>u_3>t_1>t_2>t_3>t_4>u_4>0$
\item $u_1>u_2>u_3>t_1>t_2>t_3>u_4>t_4>0$

***

\item $u_1>u_2>u_3>t_1>t_2>u_4>t_3>t_4>0$

***
***

\item $u_1>u_2>u_3>t_1>u_4>t_2>t_3>t_4>0$

***
***
***

\item $u_1>u_2>u_3>u_4>t_1>t_2>t_3>t_4>0$

\end{enumerate}

\subsection{Shuffle product of $\zeta^{(1)}(2,2) \zeta^{(1)}(2,2)$}

\begin{eqnarray*}
&&\zeta(2,2) \zeta(2,2)\\ 
& = & \sum_{\alpha_1 \in \mbox {\scriptsize{cone}}} \frac{1}{(\alpha_1)^2 (\alpha_1 + \alpha_2)^2} \sum_{\alpha_3 \in \mbox {\scriptsize{cone}}} \frac{1}{(\alpha_3)^2 (\alpha_3 + \alpha_4)^2} \\ 
& = & \sum \int_{t_1 > t_{1,2} > t_2 >t_{2,2} > t_3 > t_{3,2} > t_4 > t_{4,2} > 0} \exp (-\alpha_1t_1 - \alpha_2t_2) \exp (-\alpha_3t_3 - \alpha_4t_4) \\  
& = & \sum \frac{1}{(\alpha_1)^2 ( \alpha_1 + \alpha_2 )^2 ( \alpha_1 + \alpha_2 + \alpha_3 )^2 ( \alpha_1 + \alpha_2 + \alpha_3 + \alpha_4)^2} \\
&& +  \sum \frac{1}{(\alpha_1)^2 ( \alpha_1 + \alpha_2 )( \alpha_1 + \alpha_2 + \alpha_3 )^3 ( \alpha_1 + \alpha_2 + \alpha_3 + \alpha_4)^2} \\
&& +  \sum \frac{1}{(\alpha_1)^2 ( \alpha_1 + \alpha_2 )( \alpha_1 + \alpha_2 + \alpha_3 )^3 ( \alpha_1 + \alpha_2 + \alpha_3 + \alpha_4)^2} \\
&& +  \sum \frac{1}{(\alpha_1)^2 ( \alpha_1 + \alpha_2 )( \alpha_1 + \alpha_2 + \alpha_3 )^2 ( \alpha_1 + \alpha_2 + \alpha_3 + \alpha_4)^3} \\
&& +  \sum \frac{1}{(\alpha_1)^2 ( \alpha_1 + \alpha_2 )( \alpha_1 + \alpha_2 + \alpha_3 )^2 ( \alpha_1 + \alpha_2 + \alpha_3 + \alpha_4)^3} \\
&& +  \sum \frac{1}{(\alpha_1)^2 ( \alpha_1 + \alpha_3 )( \alpha_1 + \alpha_2 + \alpha_3 )^3 ( \alpha_1 + \alpha_2 + \alpha_3 + \alpha_4)^2} \\
&& +  \sum \frac{1}{(\alpha_1)^2 ( \alpha_1 + \alpha_3 )( \alpha_1 + \alpha_2 + \alpha_3 )^3 ( \alpha_1 + \alpha_2 + \alpha_3 + \alpha_4)^2} \\
&& +  \sum \frac{1}{(\alpha_1)^2 ( \alpha_1 + \alpha_3 )( \alpha_1 + \alpha_2 + \alpha_3 )^2 ( \alpha_1 + \alpha_2 + \alpha_3 + \alpha_4)^3} \\
&& +  \sum \frac{1}{(\alpha_1)^2 ( \alpha_1 + \alpha_3 )( \alpha_1 + \alpha_2 + \alpha_3 )^2 ( \alpha_1 + \alpha_2 + \alpha_3 + \alpha_4)^3} \\
&& +  \sum \frac{1}{(\alpha_1)^2 ( \alpha_1 + \alpha_3 )^2 ( \alpha_1 + \alpha_2 + \alpha_3 )^2 ( \alpha_1 + \alpha_2 + \alpha_3 + \alpha_4)^2} \\
&& +  \sum \frac{1}{(\alpha_1)^2 ( \alpha_1 + \alpha_3 )^2 ( \alpha_1 + \alpha_2 + \alpha_3 )( \alpha_1 + \alpha_2 + \alpha_3 + \alpha_4)^3} \\
&& +  \sum \frac{1}{(\alpha_1)^2 ( \alpha_1 + \alpha_3 )^2 ( \alpha_1 + \alpha_2 + \alpha_3 )( \alpha_1 + \alpha_2 + \alpha_3 + \alpha_4)^3} \\
&& +  \sum \frac{1}{(\alpha_1)^2 ( \alpha_1 + \alpha_3 )^2 ( \alpha_1 + \alpha_3 + \alpha_4 )( \alpha_1 + \alpha_2 + \alpha_3 + \alpha_4)^3} \\
&& +  \sum \frac{1}{(\alpha_1)^2 ( \alpha_1 + \alpha_3 )^2 ( \alpha_1 + \alpha_3 + \alpha_4 )( \alpha_1 + \alpha_2 + \alpha_3 + \alpha_4)^3} \\
&& +  \sum \frac{1}{(\alpha_1)^2 ( \alpha_1 + \alpha_3 )^2 ( \alpha_1 + \alpha_3 + \alpha_4 )^2 ( \alpha_1 + \alpha_2 + \alpha_3 + \alpha_4)^2} \\
&& +  \sum \frac{1}{(\alpha_1)( \alpha_1 + \alpha_3 )^2 ( \alpha_1 + \alpha_2 + \alpha_3)^3 (\alpha_1 + \alpha_2 + \alpha_3 + \alpha_4)^2} \\
&& +  \sum \frac{1}{(\alpha_1)( \alpha_1 + \alpha_3 )^2 ( \alpha_1 + \alpha_2 + \alpha_3)^3 (\alpha_1 + \alpha_2 + \alpha_3 + \alpha_4)^2} \\
&& +  \sum \frac{1}{(\alpha_1)( \alpha_1 + \alpha_3 )^2 ( \alpha_1 + \alpha_2 + \alpha_3)^2 (\alpha_1 + \alpha_2 + \alpha_3 + \alpha_4)^3} \\
&& +  \sum \frac{1}{(\alpha_1)( \alpha_1 + \alpha_3 )^2 ( \alpha_1 + \alpha_2 + \alpha_3)^2 (\alpha_1 + \alpha_2 + \alpha_3 + \alpha_4)^3} \\
&& +  \sum \frac{1}{(\alpha_1)( \alpha_1 + \alpha_3 )^3 ( \alpha_1 + \alpha_2 + \alpha_3)^2 (\alpha_1 + \alpha_2 + \alpha_3 + \alpha_4)^2} \\
&& +  \sum \frac{1}{(\alpha_1)( \alpha_1 + \alpha_3 )^3 ( \alpha_1 + \alpha_2 + \alpha_3)(\alpha_1 + \alpha_2 + \alpha_3 + \alpha_4)^3} \\
&& +  \sum \frac{1}{(\alpha_1)( \alpha_1 + \alpha_3 )^3 ( \alpha_1 + \alpha_2 + \alpha_3)(\alpha_1 + \alpha_2 + \alpha_3 + \alpha_4)^3} \\
&& +  \sum \frac{1}{(\alpha_1)( \alpha_1 + \alpha_3 )^3 ( \alpha_1 + \alpha_3 + \alpha_4)(\alpha_1 + \alpha_2 + \alpha_3 + \alpha_4)^3} \\
&& +  \sum \frac{1}{(\alpha_1)( \alpha_1 + \alpha_3 )^3 ( \alpha_1 + \alpha_3 + \alpha_4)(\alpha_1 + \alpha_2 + \alpha_3 + \alpha_4)^3} \\
&& +  \sum \frac{1}{(\alpha_1)( \alpha_1 + \alpha_3 )^3 ( \alpha_1 + \alpha_3 + \alpha_4)^2 (\alpha_1 + \alpha_2 + \alpha_3 + \alpha_4)^2} \\
&& +  \sum \frac{1}{(\alpha_1)( \alpha_1 + \alpha_3 )^3 ( \alpha_1 + \alpha_2 + \alpha_3)^2 (\alpha_1 + \alpha_2 + \alpha_3 + \alpha_4)^2} \\
&& +  \sum \frac{1}{(\alpha_1)( \alpha_1 + \alpha_3 )^3 ( \alpha_1 + \alpha_2 + \alpha_3)(\alpha_1 + \alpha_2 + \alpha_3 + \alpha_4)^3} \\
&& +  \sum \frac{1}{(\alpha_1)( \alpha_1 + \alpha_3 )^3 ( \alpha_1 + \alpha_2 + \alpha_3)(\alpha_1 + \alpha_2 + \alpha_3 + \alpha_4)^3} \\
&& +  \sum \frac{1}{(\alpha_1)( \alpha_1 + \alpha_3 )^3 ( \alpha_1 + \alpha_3 + \alpha_4)(\alpha_1 + \alpha_2 + \alpha_3 + \alpha_4)^3} \\
&& +  \sum \frac{1}{(\alpha_1)( \alpha_1 + \alpha_3 )^3 ( \alpha_1 + \alpha_3 + \alpha_4)(\alpha_1 + \alpha_2 + \alpha_3 + \alpha_4)^3} \\
&& +  \sum \frac{1}{(\alpha_1)( \alpha_1 + \alpha_3 )^3 ( \alpha_1 + \alpha_3 + \alpha_4)^2 (\alpha_1 + \alpha_2 + \alpha_3 + \alpha_4)^2} \\
&& +  \sum \frac{1}{(\alpha_1)( \alpha_1 + \alpha_3 )^2 ( \alpha_1 + \alpha_3 + \alpha_4)^2 (\alpha_1 + \alpha_2 + \alpha_3 + \alpha_4)^3} \\
&& +  \sum \frac{1}{(\alpha_1)( \alpha_1 + \alpha_3 )^2 ( \alpha_1 + \alpha_3 + \alpha_4)^2 (\alpha_1 + \alpha_2 + \alpha_3 + \alpha_4)^3} \\
&& +  \sum \frac{1}{(\alpha_1)( \alpha_1 + \alpha_3 )^2 ( \alpha_1 + \alpha_3 + \alpha_4)^3 (\alpha_1 + \alpha_2 + \alpha_3 + \alpha_4)^2} \\
&& +  \sum \frac{1}{(\alpha_1)( \alpha_1 + \alpha_3 )^2 ( \alpha_1 + \alpha_3 + \alpha_4)^3 (\alpha_1 + \alpha_2 + \alpha_3 + \alpha_4)^2} \\
&& +  \sum \frac{1}{(\alpha_3)( \alpha_1 + \alpha_3 )^2 ( \alpha_1 + \alpha_2 + \alpha_3)^3 (\alpha_1 + \alpha_2 + \alpha_3 + \alpha_4)^2} \\
&& +  \sum \frac{1}{(\alpha_3)( \alpha_1 + \alpha_3 )^2 ( \alpha_1 + \alpha_2 + \alpha_3)^3 (\alpha_1 + \alpha_2 + \alpha_3 + \alpha_4)^2} \\
&& +  \sum \frac{1}{(\alpha_3)( \alpha_1 + \alpha_3 )^2 ( \alpha_1 + \alpha_2 + \alpha_3)^2 (\alpha_1 + \alpha_2 + \alpha_3 + \alpha_4)^3} \\
&& +  \sum \frac{1}{(\alpha_3)( \alpha_1 + \alpha_3 )^2 ( \alpha_1 + \alpha_2 + \alpha_3)^2 (\alpha_1 + \alpha_2 + \alpha_3 + \alpha_4)^3} \\
&& +  \sum \frac{1}{(\alpha_3)( \alpha_1 + \alpha_3 )^3 ( \alpha_1 + \alpha_2 + \alpha_3)^2 (\alpha_1 + \alpha_2 + \alpha_3 + \alpha_4)^2} \\
&& +  \sum \frac{1}{(\alpha_3)( \alpha_1 + \alpha_3 )^3 ( \alpha_1 + \alpha_2 + \alpha_3)^2 (\alpha_1 + \alpha_2 + \alpha_3 + \alpha_4)^2} \\
&& +  \sum \frac{1}{(\alpha_3)( \alpha_1 + \alpha_3 )^3 ( \alpha_1 + \alpha_2 + \alpha_3)(\alpha_1 + \alpha_2 + \alpha_3 + \alpha_4)^3} \\
&& +  \sum \frac{1}{(\alpha_3)( \alpha_1 + \alpha_3 )^3 ( \alpha_1 + \alpha_3 + \alpha_4)(\alpha_1 + \alpha_2 + \alpha_3 + \alpha_4)^3} \\
&& +  \sum \frac{1}{(\alpha_3)( \alpha_1 + \alpha_3 )^3 ( \alpha_1 + \alpha_3 + \alpha_4)(\alpha_1 + \alpha_2 + \alpha_3 + \alpha_4)^3} \\
&& +  \sum \frac{1}{(\alpha_3)( \alpha_1 + \alpha_3 )^3 ( \alpha_1 + \alpha_3 + \alpha_4)^2 (\alpha_1 + \alpha_2 + \alpha_3 + \alpha_4)^2} \\
&& +  \sum \frac{1}{(\alpha_3)( \alpha_1 + \alpha_3 )^3 ( \alpha_1 + \alpha_2 + \alpha_3)^2 (\alpha_1 + \alpha_2 + \alpha_3 + \alpha_4)^2} \\
&& +  \sum \frac{1}{(\alpha_3)( \alpha_1 + \alpha_3 )^3 ( \alpha_1 + \alpha_2 + \alpha_3)(\alpha_1 + \alpha_2 + \alpha_3 + \alpha_4)^3} \\
&& +  \sum \frac{1}{(\alpha_3)( \alpha_1 + \alpha_3 )^3 ( \alpha_1 + \alpha_2 + \alpha_3)(\alpha_1 + \alpha_2 + \alpha_3 + \alpha_4)^3} \\
&& +  \sum \frac{1}{(\alpha_3)( \alpha_1 + \alpha_3 )^3 ( \alpha_1 + \alpha_2 + \alpha_3)(\alpha_1 + \alpha_2 + \alpha_3 + \alpha_4)^3} \\
&& +  \sum \frac{1}{(\alpha_3)( \alpha_1 + \alpha_3 )^3 ( \alpha_1 + \alpha_3 + \alpha_4)(\alpha_1 + \alpha_2 + \alpha_3 + \alpha_4)^3} \\
&& +  \sum \frac{1}{(\alpha_3)( \alpha_1 + \alpha_3 )^3 ( \alpha_1 + \alpha_3 + \alpha_4)(\alpha_1 + \alpha_2 + \alpha_3 + \alpha_4)^3} \\
&& +  \sum \frac{1}{(\alpha_3)( \alpha_1 + \alpha_3 )^3 ( \alpha_1 + \alpha_3 + \alpha_4)^2 (\alpha_1 + \alpha_2 + \alpha_3 + \alpha_4)^2} \\
&& +  \sum \frac{1}{(\alpha_3)( \alpha_1 + \alpha_3 )^2 ( \alpha_1 + \alpha_3 + \alpha_4)^2 (\alpha_1 + \alpha_2 + \alpha_3 + \alpha_4)^3} \\
&& +  \sum \frac{1}{(\alpha_3)( \alpha_1 + \alpha_3 )^2 ( \alpha_1 + \alpha_3 + \alpha_4)^2 (\alpha_1 + \alpha_2 + \alpha_3 + \alpha_4)^3} \\
&& +  \sum \frac{1}{(\alpha_3)( \alpha_1 + \alpha_3 )^2 ( \alpha_1 + \alpha_3 + \alpha_4)^3 (\alpha_1 + \alpha_2 + \alpha_3 + \alpha_4)^2} \\
&& +  \sum \frac{1}{(\alpha_3)( \alpha_1 + \alpha_3 )^2 ( \alpha_1 + \alpha_3 + \alpha_4)^3 (\alpha_1 + \alpha_2 + \alpha_3 + \alpha_4)^2} \\
&& +  \sum \frac{1}{(\alpha_3)^2 ( \alpha_1 + \alpha_3 )^2 ( \alpha_1 + \alpha_2 + \alpha_3)^2 (\alpha_1 + \alpha_2 + \alpha_3 + \alpha_4)^2} \\
&& +  \sum \frac{1}{(\alpha_3)^2 ( \alpha_1 + \alpha_3 )^2 ( \alpha_1 + \alpha_2 + \alpha_3)(\alpha_1 + \alpha_2 + \alpha_3 + \alpha_4)^3} \\
&& +  \sum \frac{1}{(\alpha_3)^2 ( \alpha_1 + \alpha_3 )^2 ( \alpha_1 + \alpha_2 + \alpha_3)(\alpha_1 + \alpha_2 + \alpha_3 + \alpha_4)^3} \\
&& +  \sum \frac{1}{(\alpha_3)^2 ( \alpha_1 + \alpha_3 )^2 ( \alpha_1 + \alpha_3 + \alpha_4)(\alpha_1 + \alpha_2 + \alpha_3 + \alpha_4)^3} \\
&& +  \sum \frac{1}{(\alpha_3)^2 ( \alpha_1 + \alpha_3 )^2 ( \alpha_1 + \alpha_3 + \alpha_4)(\alpha_1 + \alpha_2 + \alpha_3 + \alpha_4)^3} \\
&& +  \sum \frac{1}{(\alpha_3)^2 ( \alpha_1 + \alpha_3 )^2 ( \alpha_1 + \alpha_3 + \alpha_4)^2 (\alpha_1 + \alpha_2 + \alpha_3 + \alpha_4)^2} \\
&& +  \sum \frac{1}{(\alpha_3)^2 ( \alpha_1 + \alpha_3 )( \alpha_1 + \alpha_3 + \alpha_4)^2 (\alpha_1 + \alpha_2 + \alpha_3 + \alpha_4)^3} \\
&& +  \sum \frac{1}{(\alpha_3)^2 ( \alpha_1 + \alpha_3 )( \alpha_1 + \alpha_3 + \alpha_4)^2 (\alpha_1 + \alpha_2 + \alpha_3 + \alpha_4)^3} \\
&& +  \sum \frac{1}{(\alpha_3)^2 ( \alpha_1 + \alpha_3 )( \alpha_1 + \alpha_3 + \alpha_4)^3 (\alpha_1 + \alpha_2 + \alpha_3 + \alpha_4)^2} \\
&& +  \sum \frac{1}{(\alpha_3)^2 ( \alpha_3 + \alpha_4 )( \alpha_1 + \alpha_3 + \alpha_4)^2 (\alpha_1 + \alpha_2 + \alpha_3 + \alpha_4)^3} \\
&& +  \sum \frac{1}{(\alpha_3)^2 ( \alpha_3 + \alpha_4 )( \alpha_1 + \alpha_3 + \alpha_4)^2 (\alpha_1 + \alpha_2 + \alpha_3 + \alpha_4)^3} \\
&& +  \sum \frac{1}{(\alpha_3)^2 ( \alpha_3 + \alpha_4 )( \alpha_1 + \alpha_3 + \alpha_4)^3 (\alpha_1 + \alpha_2 + \alpha_3 + \alpha_4)^2} \\
&& +  \sum \frac{1}{(\alpha_3)^2 ( \alpha_3 + \alpha_4 )( \alpha_1 + \alpha_3 + \alpha_4)^3 (\alpha_1 + \alpha_2 + \alpha_3 + \alpha_4)^2} \\
&& +  
\sum \frac{1}{(\alpha_3)^2 ( \alpha_3 + \alpha_4 )^2 ( \alpha_1 + \alpha_3 + \alpha_4)^2 (\alpha_1 + \alpha_2 + \alpha_3 + \alpha_4)^2}= \\
& = & \zeta^{(1)}(2,2;2,2) + 2\zeta^{(1)}(2,1;3,2) + 2\zeta^{(1)}(2,1;2,3) + 2\zeta^{(23)}(2,1;3,2) \\
&& +  2\zeta^{(23)}(2,1;2,3) + \zeta^{(23)}(2,2;2,2) + 2\zeta^{(23)}(2,2;1,3) + 2\zeta^{(234)}(2,2;1,3) \\
&& +  \zeta^{(234)}(2,2;2,2) + 4\zeta^{(132)}(1,2;3,2) + 4\zeta^{(132)}(1,2;2,3) + 5\zeta^{(132)}(1,3;2,2) \\
&& +  8\zeta^{(132)}(1,3;1,3) + 7\zeta^{(1342)}(1,3;1,3) + 4\zeta^{(1342)}(1,3;2,2) + 4\zeta^{(1342)}(1,2;2,3) \\
&& +  4\zeta^{(1342)}(1,2;3,2) + \zeta^{(132)}(2,2;2,2) + 2\zeta^{(132)}(2,2;1,3) + 4\zeta^{(1342)}(2,2;1,3) \\
& &+  \zeta^{(1342)}(2,2;2,2) + 2\zeta^{(1342)}(2,1;3,2) + 2\zeta^{(13)(24)}(2,1;2,3) \\
&& +  2\zeta^{(13)(24)}(2,1;3,2) + \zeta^{(13)(24)}(2,2;2,2) .  
\end{eqnarray*}

If the cone $C$ is Galois invariant then
\begin{eqnarray*}
\zeta^{(1)}(2,2)\zeta^{(1)}(2,2)& = & 2\zeta^{(1)}(2,2;2,2)  + 4\zeta^{(23)}(2,2;2,2) \\
&& + 4\zeta^{(1)}(2,1;3,2)  + 4\zeta^{(23)}(2,1;3,2) \\
&&+4\zeta^{(1)}(2,1;2,3) + 2\zeta^{(23)}(2,1;2,3)\\
&& + 10\zeta^{(23)}(2,2;1,3) \\
&&+ 8\zeta^{(132)}(1,2;3,2) + 8\zeta^{(132)}(1,2;2,3) \\
&&+  9\zeta^{(132)}(1,3;2,2)+ 15\zeta^{(132)}(1,3;1,3). 
\end{eqnarray*}

\subsection{Shuffle product of $\zeta^{(1)}(2,2)\zeta^{(1)}(1,3)$}

\begin{eqnarray*}
&&\zeta(2,2)\zeta(1,3)=\\ 
& = & \sum_{\alpha_1 \in \mbox {\scriptsize{cone}}} \frac{1}{(\alpha_1)^2 (\alpha_1 + \alpha_2)^2} \sum_{\alpha_3 \in \mbox {\scriptsize{cone}}} \frac{1}{\alpha_3 (\alpha_3 + \alpha_4)^3} \\ 
& = & \sum \int_{t_1 > t_{2} > t_3>t_4 > 0,\,\, u_1 > u_2 > u_3 > u_4 > 0} \exp (-\alpha_1t_1 - \alpha_2t_3-\alpha_3u_1 - \alpha_4u_2)DT_8 \\ 
& = & \sum \frac{1}{(\alpha_1)^2 (\alpha_1 + \alpha_2)^2 (\alpha_1 + \alpha_2 + \alpha_3) (\alpha_1 + \alpha_2 + \alpha_3 + \alpha_4)^3} \\
&& +  \sum \frac{1}{(\alpha_1)^2 (\alpha_1 + \alpha_2) (\alpha_1 + \alpha_2 + \alpha_3)^2 (\alpha_1 + \alpha_2 + \alpha_3 + \alpha_4)^3} \\
&& +  \sum \frac{1}{(\alpha_1)^2 (\alpha_1 + \alpha_2) (\alpha_1 + \alpha_2 + \alpha_3) (\alpha_1 + \alpha_2 + \alpha_3 + \alpha_4)^4} \\
&& +  \sum \frac{1}{(\alpha_1)^2 (\alpha_1 + \alpha_2) (\alpha_1 + \alpha_2 + \alpha_3) (\alpha_1 + \alpha_2 + \alpha_3 + \alpha_4)^4} \\
&& +  \sum \frac{1}{(\alpha_1)^2 (\alpha_1 + \alpha_2) (\alpha_1 + \alpha_2 + \alpha_3) (\alpha_1 + \alpha_2 + \alpha_3 + \alpha_4)^4} \\
&& (*) +  \sum \frac{1}{(\alpha_1)^2 (\alpha_1 + \alpha_3) (\alpha_1 + \alpha_2 + \alpha_3)^2 (\alpha_1 + \alpha_2 + \alpha_3 + \alpha_4)^3} \\
&& +  \sum \frac{1}{(\alpha_1)^2 (\alpha_1 + \alpha_3) (\alpha_1 + \alpha_2 + \alpha_3)^1 (\alpha_1 + \alpha_2 + \alpha_3 + \alpha_4)^4} \\
&& +  \sum \frac{1}{(\alpha_1)^2 (\alpha_1 + \alpha_3) (\alpha_1 + \alpha_2 + \alpha_3) (\alpha_1 + \alpha_2 + \alpha_3 + \alpha_4)^4} \\
&& +  \sum \frac{1}{(\alpha_1)^2 (\alpha_1 + \alpha_3) (\alpha_1 + \alpha_2 + \alpha_3) (\alpha_1 + \alpha_2 + \alpha_3 + \alpha_4)^4} \\
&& (*) + \sum \frac{1}{(\alpha_1)^2 (\alpha_1 + \alpha_3) (\alpha_1 + \alpha_3 + \alpha_4) (\alpha_1 + \alpha_2 + \alpha_3 + \alpha_4)^4} \\
&& +  \sum \frac{1}{(\alpha_1)^2 (\alpha_1 + \alpha_3) (\alpha_1 + \alpha_3 + \alpha_4) (\alpha_1 + \alpha_2 + \alpha_3 + \alpha_4)^4} \\
&& +  \sum \frac{1}{(\alpha_1)^2 (\alpha_1 + \alpha_3) (\alpha_1 + \alpha_3 + \alpha_4) (\alpha_1 + \alpha_2 + \alpha_3 + \alpha_4)^4} \\
&& (*) +  \sum \frac{1}{(\alpha_1)^2 (\alpha_1 + \alpha_3) (\alpha_1 + \alpha_3 + \alpha_4)^2 (\alpha_1 + \alpha_2 + \alpha_3 + \alpha_4)^3} \\
&& +  \sum \frac{1}{(\alpha_1)^2 (\alpha_1 + \alpha_3) (\alpha_1 + \alpha_3 + \alpha_4)^2 (\alpha_1 + \alpha_2 + \alpha_3 + \alpha_4)^3} \\
&&(*)  +  \sum \frac{1}{(\alpha_1)^2 (\alpha_1 + \alpha_3) (\alpha_1 + \alpha_3 + \alpha_4)^3 (\alpha_1 + \alpha_2 + \alpha_3 + \alpha_4)^2} \\
&& (**) +  \sum \frac{1}{(\alpha_1) (\alpha_1 + \alpha_3)^2 (\alpha_1 + \alpha_2 + \alpha_3)^2 (\alpha_1 + \alpha_2 + \alpha_3 + \alpha_4)^3} \\
&& +  \sum \frac{1}{(\alpha_1) (\alpha_1 + \alpha_3)^2 (\alpha_1 + \alpha_2 + \alpha_3) (\alpha_1 + \alpha_2 + \alpha_3 + \alpha_4)^4} \\
&& +  \sum \frac{1}{(\alpha_1) (\alpha_1 + \alpha_3)^2 (\alpha_1 + \alpha_2 + \alpha_3) (\alpha_1 + \alpha_2 + \alpha_3 + \alpha_4)^4} \\
&& +  \sum \frac{1}{(\alpha_1) (\alpha_1 + \alpha_3)^2 (\alpha_1 + \alpha_2 + \alpha_3) (\alpha_1 + \alpha_2 + \alpha_3 + \alpha_4)^4} \\
&& (*) +  \sum \frac{1}{(\alpha_1) (\alpha_1 + \alpha_3)^2 (\alpha_1 + \alpha_3 + \alpha_4) (\alpha_1 + \alpha_2 + \alpha_3 + \alpha_4)^4} \\
&& +  \sum \frac{1}{(\alpha_1) (\alpha_1 + \alpha_3)^2 (\alpha_1 + \alpha_3 + \alpha_4) (\alpha_1 + \alpha_2 + \alpha_3 + \alpha_4)^4} \\
&& +  \sum \frac{1}{(\alpha_1) (\alpha_1 + \alpha_3)^2 (\alpha_1 + \alpha_3 + \alpha_4) (\alpha_1 + \alpha_2 + \alpha_3 + \alpha_4)^4} \\
&& (*) +  \sum \frac{1}{(\alpha_1) (\alpha_1 + \alpha_3)^2 (\alpha_1 + \alpha_3 + \alpha_4)^2 (\alpha_1 + \alpha_2 + \alpha_3 + \alpha_4)^3} \\
&& +  \sum \frac{1}{(\alpha_1) (\alpha_1 + \alpha_3)^2 (\alpha_1 + \alpha_3 + \alpha_4)^2 (\alpha_1 + \alpha_2 + \alpha_3 + \alpha_4)^3} \\
&& (*) +  \sum \frac{1}{(\alpha_1) (\alpha_1 + \alpha_3)^2 (\alpha_1 + \alpha_3 + \alpha_4)^3 (\alpha_1 + \alpha_2 + \alpha_3 + \alpha_4)^2} \\
&& (**) +  \sum \frac{1}{(\alpha_1) (\alpha_1 + \alpha_3) (\alpha_1 + \alpha_3 + \alpha_4)^2 (\alpha_1 + \alpha_2 + \alpha_3 + \alpha_4)^4} \\
&& +  \sum \frac{1}{(\alpha_1) (\alpha_1 + \alpha_3) (\alpha_1 + \alpha_3 + \alpha_4)^2 (\alpha_1 + \alpha_2 + \alpha_3 + \alpha_4)^4} \\
&& +  \sum \frac{1}{(\alpha_1) (\alpha_1 + \alpha_3)  (\alpha_1 + \alpha_3 + \alpha_4)^2 (\alpha_1 + \alpha_2 + \alpha_3 + \alpha_4)^4} \\
&& (*) +  \sum \frac{1}{(\alpha_1) (\alpha_1 + \alpha_3) (\alpha_1 + \alpha_3 + \alpha_4)^3 (\alpha_1 + \alpha_2 + \alpha_3 + \alpha_4)^3} \\
&& +  \sum \frac{1}{(\alpha_1) (\alpha_1 + \alpha_3) (\alpha_1 + \alpha_3 + \alpha_4)^3 (\alpha_1 + \alpha_2 + \alpha_3 + \alpha_4)^3} \\
&& (*) +  \sum \frac{1}{(\alpha_1) (\alpha_1 + \alpha_3) (\alpha_1 + \alpha_3 + \alpha_4)^4 (\alpha_1 + \alpha_2 + \alpha_3 + \alpha_4)^2} \\
&& (**) +  \sum \frac{1}{(\alpha_1) (\alpha_1 + \alpha_3) (\alpha_1 + \alpha_3 + \alpha_4)^3 (\alpha_1 + \alpha_2 + \alpha_3 + \alpha_4)^3} \\
&& +  \sum \frac{1}{(\alpha_1) (\alpha_1 + \alpha_3) (\alpha_1 + \alpha_3 + \alpha_4)^3 (\alpha_1 + \alpha_2 + \alpha_3 + \alpha_4)^3} \\
&& (*) +  \sum \frac{1}{(\alpha_1) (\alpha_1 + \alpha_3) (\alpha_1 + \alpha_3 + \alpha_4)^4 (\alpha_1 + \alpha_2 + \alpha_3 + \alpha_4)^2} \\
&& (**) +  \sum \frac{1}{(\alpha_1) (\alpha_1 + \alpha_3) (\alpha_1 + \alpha_3 + \alpha_4)^4 (\alpha_1 + \alpha_2 + \alpha_3 + \alpha_4)^2} \\
&& (***)+  \sum \frac{1}{(\alpha_3) (\alpha_1 + \alpha_3)^2 (\alpha_1 + \alpha_2 + \alpha_3)^2 (\alpha_1 + \alpha_2 + \alpha_3 + \alpha_4)^3} \\
&& +  \sum \frac{1}{(\alpha_3) (\alpha_1 + \alpha_3)^2 (\alpha_1 + \alpha_2 + \alpha_3) (\alpha_1 + \alpha_2 + \alpha_3 + \alpha_4)^4} \\
&& +  \sum \frac{1}{(\alpha_3) (\alpha_1 + \alpha_3)^2 (\alpha_1 + \alpha_2 + \alpha_3) (\alpha_1 + \alpha_2 + \alpha_3 + \alpha_4)^4} \\
&& +  \sum \frac{1}{(\alpha_3) (\alpha_1 + \alpha_3)^2 (\alpha_1 + \alpha_2 + \alpha_3) (\alpha_1 + \alpha_2 + \alpha_3 + \alpha_4)^4} \\
&& (*) +  \sum \frac{1}{(\alpha_3) (\alpha_1 + \alpha_3)^2 (\alpha_1 + \alpha_3 + \alpha_4) (\alpha_1 + \alpha_2 + \alpha_3 + \alpha_4)^4} \\
&& +  \sum \frac{1}{(\alpha_3) (\alpha_1 + \alpha_3)^2 (\alpha_1 + \alpha_3 + \alpha_4) (\alpha_1 + \alpha_2 + \alpha_3 + \alpha_4)^4} \\
&& +  \sum \frac{1}{(\alpha_3) (\alpha_1 + \alpha_3)^2 (\alpha_1 + \alpha_3 + \alpha_4) (\alpha_1 + \alpha_2 + \alpha_3 + \alpha_4)^4} \\
&& (*)+  \sum \frac{1}{(\alpha_3) (\alpha_1 + \alpha_3)^2 (\alpha_1 + \alpha_3 + \alpha_4)^2 (\alpha_1 + \alpha_2 + \alpha_3 + \alpha_4)^3} \\
&& +  \sum \frac{1}{(\alpha_3) (\alpha_1 + \alpha_3)^2 (\alpha_1 + \alpha_3 + \alpha_4)^2 (\alpha_1 + \alpha_2 + \alpha_3 + \alpha_4)^3} \\
&& (*)+  \sum \frac{1}{(\alpha_3) (\alpha_1 + \alpha_3)^2 (\alpha_1 + \alpha_3 + \alpha_4)^3 (\alpha_1 + \alpha_2 + \alpha_3 + \alpha_4)^2} \\
&& (**)+  \sum \frac{1}{(\alpha_3) (\alpha_1 + \alpha_3) (\alpha_1 + \alpha_3 + \alpha_4)^2 (\alpha_1 + \alpha_2 + \alpha_3 + \alpha_4)^4} \\
&& +  \sum \frac{1}{(\alpha_3) (\alpha_1 + \alpha_3) (\alpha_1 + \alpha_3 + \alpha_4)^2 (\alpha_1 + \alpha_2 + \alpha_3 + \alpha_4)^4} \\
&& +  \sum \frac{1}{(\alpha_3) (\alpha_1 + \alpha_3) (\alpha_1 + \alpha_3 + \alpha_4)^2 (\alpha_1 + \alpha_2 + \alpha_3 + \alpha_4)^4} \\
&& (*)+  \sum \frac{1}{(\alpha_3) (\alpha_1 + \alpha_3) (\alpha_1 + \alpha_3 + \alpha_4)^3 (\alpha_1 + \alpha_2 + \alpha_3 + \alpha_4)^3} \\
&& +  \sum \frac{1}{(\alpha_3) (\alpha_1 + \alpha_3) (\alpha_1 + \alpha_3 + \alpha_4)^3 (\alpha_1 + \alpha_2 + \alpha_3 + \alpha_4)^3} \\
&& (*)+  \sum \frac{1}{(\alpha_3) (\alpha_1 + \alpha_3) (\alpha_1 + \alpha_3 + \alpha_4)^4 (\alpha_1 + \alpha_2 + \alpha_3 + \alpha_4)^2} \\
&& (**)+  \sum \frac{1}{(\alpha_3) (\alpha_1 + \alpha_3) (\alpha_1 + \alpha_3 + \alpha_4)^3 (\alpha_1 + \alpha_2 + \alpha_3 + \alpha_4)^3} \\
&& +  \sum \frac{1}{(\alpha_3) (\alpha_1 + \alpha_3) (\alpha_1 + \alpha_3 + \alpha_4)^3 (\alpha_1 + \alpha_2 + \alpha_3 + \alpha_4)^3} \\
&& (*)+  \sum \frac{1}{(\alpha_3) (\alpha_1 + \alpha_3) (\alpha_1 + \alpha_3 + \alpha_4)^2 (\alpha_1 + \alpha_2 + \alpha_3 + \alpha_4)^4} \\
&& (**)+  \sum \frac{1}{(\alpha_3) (\alpha_1 + \alpha_3) (\alpha_1 + \alpha_3 + \alpha_4)^4 (\alpha_1 + \alpha_2 + \alpha_3 + \alpha_4)^2} \\
&& (***)+  \sum \frac{1}{(\alpha_3) (\alpha_3 + \alpha_4) (\alpha_1 + \alpha_3 + \alpha_4)^2 (\alpha_1 + \alpha_2 + \alpha_3 + \alpha_4)^4} \\
&& +  \sum \frac{1}{(\alpha_3) (\alpha_3 + \alpha_4) (\alpha_1 + \alpha_3 + \alpha_4)^2 (\alpha_1 + \alpha_2 + \alpha_3 + \alpha_4)^4} \\
&& +  \sum \frac{1}{(\alpha_3) (\alpha_3 + \alpha_4) (\alpha_1 + \alpha_3 + \alpha_4)^2 (\alpha_1 + \alpha_2 + \alpha_3 + \alpha_4)^4} \\
&& (*)+  \sum \frac{1}{(\alpha_3) (\alpha_3 + \alpha_4) (\alpha_1 + \alpha_3 + \alpha_4)^3 (\alpha_1 + \alpha_2 + \alpha_3 + \alpha_4)^3} \\
&& +  \sum \frac{1}{(\alpha_3) (\alpha_3 + \alpha_4) (\alpha_1 + \alpha_3 + \alpha_4)^3 (\alpha_1 + \alpha_2 + \alpha_3 + \alpha_4)^3} \\
&& (*)+  \sum \frac{1}{(\alpha_3) (\alpha_3 + \alpha_4) (\alpha_1 + \alpha_3 + \alpha_4)^4 (\alpha_1 + \alpha_2 + \alpha_3 + \alpha_4)^2} \\
&& (**)+  \sum \frac{1}{(\alpha_3) (\alpha_3 + \alpha_4) (\alpha_1 + \alpha_3 + \alpha_4)^2 (\alpha_1 + \alpha_2 + \alpha_3 + \alpha_4)^4} \\
&& +  \sum \frac{1}{(\alpha_3) (\alpha_3 + \alpha_4) (\alpha_1 + \alpha_3 + \alpha_4)^2 (\alpha_1 + \alpha_2 + \alpha_3 + \alpha_4)^4} \\
&& (*)+  \sum \frac{1}{(\alpha_3) (\alpha_3 + \alpha_4) (\alpha_1 + \alpha_3 + \alpha_4)^2 (\alpha_1 + \alpha_2 + \alpha_3 + \alpha_4)^4} \\
&& (**)+  \sum \frac{1}{(\alpha_3) (\alpha_3 + \alpha_4) (\alpha_1 + \alpha_3 + \alpha_4)^3(\alpha_1 + \alpha_2 + \alpha_3 + \alpha_4)^3} \\
&& (***)+  \sum \frac{1}{(\alpha_3) (\alpha_3 + \alpha_4)^2 (\alpha_1 + \alpha_3 + \alpha_4)^2 (\alpha_1 + \alpha_2 + \alpha_3 + \alpha_4)^3} \\
&& +  \sum \frac{1}{(\alpha_3) (\alpha_3 + \alpha_4)^2 (\alpha_1 + \alpha_3 + \alpha_4)^2 (\alpha_1 + \alpha_2 + \alpha_3 + \alpha_4)^3} \\
&& (*)+  \sum \frac{1}{(\alpha_3) (\alpha_3 + \alpha_4)^2 (\alpha_1 + \alpha_3 + \alpha_4)^3 (\alpha_1 + \alpha_2 + \alpha_3 + \alpha_4)^2} \\
&& (**)+  \sum \frac{1}{(\alpha_3) (\alpha_3 + \alpha_4)^2 (\alpha_1 + \alpha_3 + \alpha_4)^3 (\alpha_1 + \alpha_2 + \alpha_3 + \alpha_4)^2} \\
&& (***)+  \sum \frac{1}{(\alpha_3) (\alpha_3 + \alpha_4)^3 (\alpha_1 + \alpha_3 + \alpha_4)^2 (\alpha_1 + \alpha_2 + \alpha_3 + \alpha_4)^2} =\\
& = & \zeta^{(1)}(2,2;1,3) + \zeta^{(1)}(2,1;2,3) +3\zeta^{(1)}(2,1;1,4)+ \zeta^{(23)}(2,1;2,3) \\
&& +  3\zeta^{(23)}(2,1;1,4) + 3\zeta^{(234)}(2,1;1,4) + 2\zeta^{(234)}(2,1;2,3) + \zeta^{(234)}(2,1;3,2) \\
&& +  \zeta^{(23)}(1,2;2,3) + 3\zeta^{(23)}(1,2;1,4) + 3\zeta^{(234)}(1,2;1,4) + 2\zeta^{(234)}(1,2;2,3) \\
&& +  \zeta^{(234)}(1,2;3,2) + 3\zeta^{(234)}(1,1;2,4) + 2\zeta^{(234)}(1,1;3,3) + \zeta^{(234)}(1,1;4,2) \\
&& +2\zeta^{(234)}(1,1,3,3) +2\zeta^{(234)}(1,1;4,2) + \zeta^{(132)}(1,2;2,3) + 3\zeta^{(132)}(1,2;1,4) \\
&&+ 3\zeta^{(1342)}(1,2;1,4) + 2\zeta^{(1342)}(1,2;2,3) + \zeta^{(1342)}(1,2;3,2)\\
&& + 3\zeta^{(1342)}(1,1;2,4) + 2\zeta^{(1342)}(1,1;3,3) +  \zeta^{(1342)}(1,1;4,2) + 2\zeta^{(1342)}(1,1;3,3)\\
&& + \zeta^{(1342)}(1,1,2,4) +  \zeta^{(1342)}(1,1,4,2)  + 3\zeta^{(13)(24)}(1,1;2,4) \\
&& +  2\zeta^{(13)(24)}(1,1;3,3) + \zeta^{(13)(24)}(1,1;4,2) +3\zeta^{(13)(24)}(1,1,2,4) 
+ \zeta^{(13)(24)}(1,1;3,3) \\
&& +  2\zeta^{(13)(24)}(1,2;2,3) + 2\zeta^{(13)(24)}(1,2;3,2) + \zeta^{(13)(24)}(1,3;2,2). 
\end{eqnarray*}

If the cone is Galois invariant then
\begin{eqnarray*}
&&\zeta^{(1)}(2,2)\zeta^{(1)}(1,3)=\\
& = & \zeta^{(1)}(2,2;1,3) + \zeta^{(1)}(2,1;2,3) + 3\zeta^{(23)}(2,1;2,3) +\\
&&+3\zeta^{(1)}(2,1;1,4)  +6\zeta^{(23)}(2,1;1,4) + \zeta^{(23)}(2,1;3,2) \\
&&   +2\zeta^{(1)}(1,2;2,3) + 6\zeta^{(23)}(1,2,2,3) + 12\zeta^{(23)}(1,2;1,4)  \\
&& +2\zeta^{(1)}(1,2,3,2) +2\zeta^{(23)}(1,2,3,2)\\
&& + 6\zeta^{(1)}(1,1,2,4) + 7\zeta^{(23)}(1,2,1,4) + 3\zeta^{(1)}(1,1,3,3) + 8\zeta^{(23)}(1,1,3,3)\\
&& + \zeta^{(1)}(1,1,4,2) + 5\zeta^{(23)}(1,1;3,3)+\zeta^{(1)}(1,3,2,2). 
\end{eqnarray*}
\subsection{Shuffle product of $\zeta^{(1)}(1,3)\zeta^{(1)}(1,3)$}

\begin{eqnarray*}
&&\zeta(1,3) \zeta(1,3)=\\
 & = & \sum_{\alpha_1 \in \mbox {\scriptsize{cone}}} \frac{1}{\alpha_1 (\alpha_1 + \alpha_2)^3} \sum_{\alpha_3 \in \mbox {\scriptsize{cone}}} \frac{1}{\alpha_3 (\alpha_3 + \alpha_4)^3} \\
& = & \sum \int_{t_1 > t_2 > t_{2,2} >t_{2,3} > 0,\,\,t_3 > t_4 > t_{4,2} > t_{4,3} > 0} \exp (-\alpha_1t_1 - \alpha_2t_2-\alpha_3t_3 - \alpha_4t_4)DT_8 \\    
& = & \zeta^{1}(1,3;1,3) + \zeta^{1}(1,2;2,3) + 3\zeta^{1}(1,2;1,4) + \zeta^{1}(1,1;3,3) + 3\zeta^{1}(1,1;2,4) \\
& + &  6\zeta^{1}(1,1;1,5) + \zeta^{(23)}(1,1;3,3) + 3\zeta^{(23)}(1,1;2,4) + 6\zeta^{(23)}(1,1;1,5) + 6\zeta^{(234)}(1,1;1,5) \\
& + & 3\zeta^{(234)}(1,1;2,4) + \zeta^{(234)}(1,1;3,3) + \zeta^{(132)}(1,1;3,3) + 3\zeta^{(132)}(1,1;2,4) + 6\zeta^{(132)}(1,1;1,5)\\
& + & 6\zeta^{(1342)}(1,1;1,5) + 3\zeta^{(1342)}(1,1;2,4) + \zeta^{(1342)}(1,1;3,3) + 6\zeta^{(13)(24)}(1,1;1,5) \\
& + & 3\zeta^{(13)(24)}(1,1;2,4) + \zeta^{(13)(24)}(1,1;3,3) + 3\zeta^{(13)(24)}(1,2;1,4) + \zeta^{(13)(24)}(1,2;2,3) \\
& + & \zeta^{(13)(24)}(1,3;1,3) .
\end{eqnarray*}

If the cone $C$ is Galois invariant then

\begin{eqnarray*}
\zeta(1,3) \zeta(1,3)& = & 2\zeta^{(1)}(1,3;1,3) + 2\zeta^{(1)}(1,2;2,3) + 6\zeta^{(1)}(1,2;1,4)\\
&& + 2\zeta^{1}(1,1;3,3)+ 4\zeta^{(23)}(1,1;3,3)\\
&&+ 6\zeta^{1}(1,1;2,4) +  12\zeta^{(23)}(1,1;2,4)  \\
&& +   12\zeta^{1}(1,1;1,5) + 24\zeta^{(23)}(1,1;1,5).
 \end{eqnarray*}

\subsection{Formula for the shuffle product of $\zeta_{K,C}(2)\zeta_{K,C}(2)$}
Using the self shuffle of  $\zeta_{K,C}(2)$  together with the simplified formulas for the previous four subsections, assuming Galois invariance of the cone, we obtain the following Theorem. 

\begin{theorem}
The shuffle product of $\zeta_{K,C}(2)$ times  itself in terms of refinements of MDZVs is given by
\begin{eqnarray*}
&&\zeta_{K,C}(2)\zeta_{K,C}(2)=\\
&=&
8\zeta^{(1)}(2,2,2,2)+
16\zeta^{(1)}(2,1,3,2)+
32\zeta^{(1)}(2,1,2,3)+
8\zeta^{(1)}(2,2,1,3)\\
&&+
24\zeta^{(1)}(2,1,1,4)+
48\zeta^{(1)}(1,2,2,3)+
16\zeta^{(1)}(1,2,3,2)\\
&&+144\zeta^{(1)}(1,1,2,4)+
56\zeta^{(1)}(1,1,3,3)+
8\zeta^{(1)}(1,1,4,2)+
8\zeta^{(1)}(1,3,2,2)\\
&&+
32\zeta^{(1)}(1,3,1,3)+
8\zeta^{(1)}(1,2,1,4)+
8\zeta^{(1)}(1,1,15)\\
&&+
16\zeta^{(23)}(2,2,2,2)+
24\zeta^{(23)}(2,1,3,2)+
32\zeta^{(23)}(2,1,2,3)+
40\zeta^{(23)}(2,2,1,3)\\
&&+
48\zeta^{(23)}(1,2,3,2)+
80\zeta^{(23)}(1,2,2,3)+
36\zeta^{(23)}(1,3,2,2)+
60\zeta^{(23)}(1,3,1,3)\\
&&+
48\zeta^{(23)}(2,1,1,4)+
96\zeta^{(23)}(1,2,1,4)+
248\zeta^{(23)}(1,1,2,4)\\
&&+
128\zeta^{(23)}(1,1,3,3)+
40\zeta^{(23)}(1,1,4,2)+
384\zeta^{(23)}(1,1,1,5).
\end{eqnarray*}

\end{theorem}

\section{Shuffle product $\zeta_{K,C}(2)\zeta_{K,C}(3)$}

\subsection{Preliminary}
In this Section we consider the shuffle of the product $\zeta_{K,C}(2)\zeta_{K,C}(3)$. In order to obtain the answer in the form of the new type of MDZVs, we proceed as follows. First, we perform a self shuffle of $\zeta_{K,C}(2)$ and of $\zeta_{K,C}(3)$.

Earlier, in Propositions \label{prop self-sh of zeta(2)} and \label{prop self-sh of zeta(3)}
we computed the selfie of Dedekind zeta $\zeta_{K,C}(2)$, which is
\[\zeta_{K,C}(2)=2\zeta^{(1)}(2,2)+4\zeta^{(1)}(1,3).\]
and the selfie of the Dedekind zeta $\zeta_{K,C}(3)$, which is

\begin{equation}
\nonumber
\zeta_{K,C} (3) = 2\zeta^{(1)}(3,3) + 6\zeta^{(1)}(2,4) + 12\zeta^{(1)}(1,5).
\end{equation}
Then we have
\[\zeta_{K,C}(2)\zeta_{K,C} (3)=\left(2\zeta^{(1)}(2,2)+4\zeta^{(1)}(1,3)\right)\left(2\zeta^{(1)}(3,3) + 6\zeta^{(1)}(2,4) + 12\zeta^{(1)}(1,5)\right).\]

In order to compute the shuffle product $\zeta_{K,C}(2)\zeta_{K,C}(3)$ in terms of the new MDZVs, we have to perform the shuffles the six shuffles:
\begin{enumerate}
\item
$\zeta^{(1)}(2,2)\zeta^{(1)}(3,3)$, 

\item
$\zeta^{(1)}(2,2)\zeta^{(1)}(2,4)$

\item
$\zeta^{(1)}(2,2)\zeta^{(1)}(1,5)$,

\item
$\zeta^{(1)}(1,3)\zeta^{(1)}(3,3)$, 

\item
$\zeta^{(1)}(1,3)\zeta^{(1)}(2,4)$

\item
$\zeta^{(1)}(1,3)\zeta^{(1)}(1,5)$.
\end{enumerate}

For the $\zeta^{(1)}(2,2)$
we use
\begin{eqnarray*}
&&\zeta^{(1)}(2,2)=\\
&&=\int_{t_1>t_2>t_3>t_4>0}\exp(-\alpha_1t_1-\alpha_2t_3)dt_1dt_2dt_3dt_4
\end{eqnarray*}

For the $\zeta^{(1)}(1,3)$
we use
\begin{eqnarray*}
&&\zeta^{(1)}(1,3)=\\
&&=\int_{t_1>t_2>t_3>t_4>0}\exp(-\alpha_1t_1-\alpha_2t_2)dt_1dt_2dt_3dt_4
\end{eqnarray*}

For the $\zeta^{(1)}(3,3)$
we use
\begin{eqnarray*}
&&\zeta^{(1)}(3,3)=\\
&&=\int_{u_1>u_2>u_3>u_4>u_5>u_6>0}\exp(-\beta_1u_1-\beta_2u_4)dt_1dt_2dt_3dt_4
\end{eqnarray*}

For the $\zeta^{(1)}(2,4)$
we use
\begin{eqnarray*}
&&\zeta^{(1)}(2,4)=\\
&&=\int_{u_1>u_2>u_3>u_4>u_5>u_6>0}\exp(-\beta_1u_1-\beta_2u_3)dt_1dt_2dt_3dt_4
\end{eqnarray*}

For the $\zeta^{(1)}(1,5)$
we use
\begin{eqnarray*}
&&\zeta^{(1)}(1,5)=\\
&&=\int_{u_1>u_2>u_3>u_4>u_5>u_6>0}\exp(-\beta_1u_1-\beta_2u_2)dt_1dt_2dt_3dt_4
\end{eqnarray*}

The following Subsections deal with the above list of six shuffles. Before we can proceed with this, it is very useful to have the list of all needed shuffles of chains of inequalities.
For that purpose the next subsection gives a list of all such inequalities.

\subsection{Needed inequalities for the shuffles}
In this subsection we give a list of all possible shuffle of the two inequalities
\[t_1>t_2>t_3>t_4>0\]
and
\[u_1>u_2>u_3>u_4>u_5>u_6>0.\]

It is possible to list such shuffles algorithmically, however that task is left for other people to consider.

\begin{enumerate}
\item $t_1>t_2>t_3>t_4>u_1>u_2>u_3>u_4>u_5>u_6>0$
\item $t_1>t_2>t_3>u_1>t_4>u_2>u_3>u_4>u_5>u_6>0$
\item $t_1>t_2>t_3>u_1>u_2>t_4>u_3>u_4>u_5>u_6>0$
\item $t_1>t_2>t_3>u_1>u_2>u_3>t_4>u_4>u_5>u_6>0$
\item $t_1>t_2>t_3>u_1>u_2>u_3>u_4>t_4>u_5>u_6>0$
\item $t_1>t_2>t_3>u_1>u_2>u_3>u_4>u_5>t_4>u_6>0$
\item $t_1>t_2>t_3>u_1>u_2>u_3>u_4>u_5>u_6>t_4>0$

***

\item $t_1>t_2>u_1>t_3>t_4>u_2>u_3>u_4>u_5>u_6>0$
\item $t_1>t_2>u_1>t_3>u_2>t_4>u_3>u_4>u_5>u_6>0$
\item $t_1>t_2>u_1>t_3>u_2>u_3>t_4>u_4>u_5>u_6>0$
\item $t_1>t_2>u_1>t_3>u_2>u_3>u_4>t_4>u_5>u_6>0$
\item $t_1>t_2>u_1>t_3>u_2>u_3>u_4>u_5>t_4>u_6>0$
\item $t_1>t_2>u_1>t_3>u_2>u_3>u_4>u_5>u_6>t_4>0$

***

\item $t_1>t_2>u_1>u_2>t_3>t_4>u_3>u_4>u_5>u_6>0$
\item $t_1>t_2>u_1>u_2>t_3>u_3>t_4>u_4>u_5>u_6>0$
\item $t_1>t_2>u_1>u_2>t_3>u_3>u_4>t_4>u_5>u_6>0$
\item $t_1>t_2>u_1>u_2>t_3>u_3>u_4>u_5>t_4>u_6>0$
\item $t_1>t_2>u_1>u_2>t_3>u_3>u_4>u_5>u_6>t_4>0$

***

\item $t_1>t_2>u_1>u_2>u_3>t_3>t_4>u_4>u_5>u_6>0$
\item $t_1>t_2>u_1>u_2>u_3>t_3>u_4>t_4>u_5>u_6>0$
\item $t_1>t_2>u_1>u_2>u_3>t_3>u_4>u_5>t_4>u_6>0$
\item $t_1>t_2>u_1>u_2>u_3>t_3>u_4>u_5>u_6>t_4>0$

***

\item $t_1>t_2>u_1>u_2>u_3>u_4>t_3>t_4>u_5>u_6>0$
\item $t_1>t_2>u_1>u_2>u_3>u_4>t_3>u_5>t_4>u_6>0$
\item $t_1>t_2>u_1>u_2>u_3>u_4>t_3>u_5>u_6>t_4>0$

***

\item $t_1>t_2>u_1>u_2>u_3>u_4>u_5>t_3>t_4>u_6>0$
\item $t_1>t_2>u_1>u_2>u_3>u_4>u_5>t_3>u_6>t_4>0$

*** 

\item $t_1>t_2>u_1>u_2>u_3>u_4>u_5>u_6>t_3>t_4>0$

***
***

\item $t_1>u_1>t_2>t_3>t_4>u_2>u_3>u_4>u_5>u_6>0$
\item $t_1>u_1>t_2>t_3>u_2>t_4>u_3>u_4>u_5>u_6>0$
\item $t_1>u_1>t_2>t_3>u_2>u_3>t_4>u_4>u_5>u_6>0$
\item $t_1>u_1>t_2>t_3>u_2>u_3>u_4>t_4>u_5>u_6>0$
\item $t_1>u_1>t_2>t_3>u_2>u_3>u_4>u_5>t_4>u_6>0$
\item $t_1>u_1>t_2>t_3>u_2>u_3>u_4>u_5>u_6>t_4>0$

***

\item $t_1>u_1>t_2>u_2>t_3>t_4>u_3>u_4>u_5>u_6>0$
\item $t_1>u_1>t_2>u_2>t_3>u_3>t_4>u_4>u_5>u_6>0$
\item $t_1>u_1>t_2>u_2>t_3>u_3>u_4>t_4>u_5>u_6>0$
\item $t_1>u_1>t_2>u_2>t_3>u_3>u_4>u_5>t_4>u_6>0$
\item $t_1>u_1>t_2>u_2>t_3>u_3>u_4>u_5>u_6>t_4>0$

***

\item $t_1>u_1>t_2>u_2>u_3>t_3>t_4>u_4>u_5>u_6>0$
\item $t_1>u_1>t_2>u_2>u_3>t_3>u_4>t_4>u_5>u_6>0$
\item $t_1>u_1>t_2>u_2>u_3>t_3>u_4>u_5>t_4>u_6>0$
\item $t_1>u_1>t_2>u_2>u_3>t_3>u_4>u_5>u_6>t_4>0$

***

\item $t_1>u_1>t_2>u_2>u_3>u_4>t_3>t_4>u_5>u_6>0$
\item $t_1>u_1>t_2>u_2>u_3>u_4>t_3>u_5>t_4>u_6>0$
\item $t_1>u_1>t_2>u_2>u_3>u_4>t_3>u_5>u_6>t_4>0$

***

\item $t_1>u_1>t_2>u_2>u_3>u_4>u_5>t_3>t_4>u_6>0$
\item $t_1>u_1>t_2>u_2>u_3>u_4>u_5>t_3>u_6>t_4>0$

***

\item $t_1>u_1>t_2>u_2>u_3>u_4>u_5>u_6>t_3>t_4>0$

***
***

\item $t_1>u_1>u_2>t_2>t_3>t_4>u_3>u_4>u_5>u_6>0$
\item $t_1>u_1>u_2>t_2>t_3>u_3>t_4>u_4>u_5>u_6>0$
\item $t_1>u_1>u_2>t_2>t_3>u_3>u_4>t_4>u_5>u_6>0$
\item $t_1>u_1>u_2>t_2>t_3>u_3>u_4>u_5>t_4>u_6>0$
\item $t_1>u_1>u_2>t_2>t_3>u_3>u_4>u_5>u_6>t_4>0$

***

\item $t_1>u_1>u_2>t_2>u_3>t_3>t_4>u_4>u_5>u_6>0$
\item $t_1>u_1>u_2>t_2>u_3>t_3>u_4>t_4>u_5>u_6>0$
\item $t_1>u_1>u_2>t_2>u_3>t_3>u_4>u_5>t_4>u_6>0$
\item $t_1>u_1>u_2>t_2>u_3>t_3>u_4>u_5>u_6>t_4>0$

***

\item $t_1>u_1>u_2>t_2>u_3>u_4>t_3>t_4>u_5>u_6>0$
\item $t_1>u_1>u_2>t_2>u_3>u_4>t_3>u_5>t_4>u_6>0$
\item $t_1>u_1>u_2>t_2>u_3>u_4>t_3>u_5>u_6>t_4>0$

***

\item $t_1>u_1>u_2>t_2>u_3>u_4>u_5>t_3>t_4>u_6>0$
\item $t_1>u_1>u_2>t_2>u_3>u_4>u_5>t_3>u_6>t_4>0$

***

\item $t_1>u_1>u_2>t_2>u_3>u_4>u_5>u_6>t_3>t_4>0$

***
***

\item $t_1>u_1>u_2>u_3>t_2>t_3>t_4>u_4>u_5>u_6>0$
\item $t_1>u_1>u_2>u_3>t_2>t_3>u_4>t_4>u_5>u_6>0$
\item $t_1>u_1>u_2>u_3>t_2>t_3>u_4>u_5>t_4>u_6>0$
\item $t_1>u_1>u_2>u_3>t_2>t_3>u_4>u_5>u_6>t_4>0$

***

\item $t_1>u_1>u_2>u_3>t_2>u_4>t_3>t_4>u_5>u_6>0$
\item $t_1>u_1>u_2>u_3>t_2>u_4>t_3>u_5>t_4>u_6>0$
\item $t_1>u_1>u_2>u_3>t_2>u_4>t_3>u_5>u_6>t_4>0$

***

\item $t_1>u_1>u_2>u_3>t_2>u_4>u_5>t_3>t_4>u_6>0$
\item $t_1>u_1>u_2>u_3>t_2>u_4>u_5>t_3>u_6>t_4>0$

***

\item $t_1>u_1>u_2>u_3>t_2>u_4>u_5>u_6>t_3>t_4>0$

***
***

\item $t_1>u_1>u_2>u_3>u_4>t_2>t_3>t_4>u_5>u_6>0$
\item $t_1>u_1>u_2>u_3>u_4>t_2>t_3>u_5>t_4>u_6>0$
\item $t_1>u_1>u_2>u_3>u_4>t_2>t_3>u_5>u_6>t_4>0$

***

\item $t_1>u_1>u_2>u_3>u_4>t_2>u_5>t_3>t_4>u_6>0$
\item $t_1>u_1>u_2>u_3>u_4>t_2>u_5>t_3>u_6>t_4>0$

***

\item $t_1>u_1>u_2>u_3>u_4>t_2>u_5>u_6>t_3>t_4>0$

***
***

\item $t_1>u_1>u_2>u_3>u_4>u_5>t_2>t_3>t_4>u_6>0$
\item $t_1>u_1>u_2>u_3>u_4>u_5>t_2>t_3>u_6>t_4>0$

***

\item $t_1>u_1>u_2>u_3>u_4>u_5>t_2>u_6>t_3>t_4>0$

***
***

\item $t_1>u_1>u_2>u_3>u_4>u_5>u_6>t_2>t_3>t_4>0$

***
***
***

\item $u_1>t_1>t_2>t_3>t_4>u_2>u_3>u_4>u_5>u_6>0$
\item $u_1>t_1>t_2>t_3>u_2>t_4>u_3>u_4>u_5>u_6>0$
\item $u_1>t_1>t_2>t_3>u_2>u_3>t_4>u_4>u_5>u_6>0$
\item $u_1>t_1>t_2>t_3>u_2>u_3>u_4>t_4>u_5>u_6>0$
\item $u_1>t_1>t_2>t_3>u_2>u_3>u_4>u_5>t_4>u_6>0$
\item $u_1>t_1>t_2>t_3>u_2>u_3>u_4>u_5>u_6>t_4>0$

***

\item $u_1>t_1>t_2>u_2>t_3>t_4>u_3>u_4>u_5>u_6>0$
\item $u_1>t_1>t_2>u_2>t_3>u_3>t_4>u_4>u_5>u_6>0$
\item $u_1>t_1>t_2>u_2>t_3>u_3>u_4>t_4>u_5>u_6>0$
\item $u_1>t_1>t_2>u_2>t_3>u_3>u_4>u_5>t_4>u_6>0$
\item $u_1>t_1>t_2>u_2>t_3>u_3>u_4>u_5>u_6>t_4>0$

***

\item $u_1>t_1>t_2>u_2>u_3>t_3>t_4>u_4>u_5>u_6>0$
\item $u_1>t_1>t_2>u_2>u_3>t_3>u_4>t_4>u_5>u_6>0$
\item $u_1>t_1>t_2>u_2>u_3>t_3>u_4>u_5>t_4>u_6>0$
\item $u_1>t_1>t_2>u_2>u_3>t_3>u_4>u_5>u_6>t_4>0$

***

\item $u_1>t_1>t_2>u_2>u_3>u_4>t_3>t_4>u_5>u_6>0$
\item $u_1>t_1>t_2>u_2>u_3>u_4>t_3>u_5>t_4>u_6>0$
\item $u_1>t_1>t_2>u_2>u_3>u_4>t_3>u_5>u_6>t_4>0$

***

\item $u_1>t_1>t_2>u_2>u_3>u_4>u_5>t_3>t_4>u_6>0$
\item $u_1>t_1>t_2>u_2>u_3>u_4>u_5>t_3>u_6>t_4>0$

***

\item $u_1>t_1>t_2>u_2>u_3>u_4>u_5>u_6>t_3>t_4>0$

***
***

\item $u_1>t_1>u_2>t_2>t_3>t_4>u_3>u_4>u_5>u_6>0$
\item $u_1>t_1>u_2>t_2>t_3>u_3>t_4>u_4>u_5>u_6>0$
\item $u_1>t_1>u_2>t_2>t_3>u_3>u_4>t_4>u_5>u_6>0$
\item $u_1>t_1>u_2>t_2>t_3>u_3>u_4>u_5>t_4>u_6>0$
\item $u_1>t_1>u_2>t_2>t_3>u_3>u_4>u_5>u_6>t_4>0$

***

\item $u_1>t_1>u_2>t_2>u_3>t_3>t_4>u_4>u_5>u_6>0$
\item $u_1>t_1>u_2>t_2>u_3>t_3>u_4>t_4>u_5>u_6>0$
\item $u_1>t_1>u_2>t_2>u_3>t_3>u_4>u_5>t_4>u_6>0$
\item $u_1>t_1>u_2>t_2>u_3>t_3>u_4>u_5>u_6>t_4>0$

***

\item $u_1>t_1>u_2>t_2>u_3>u_4>t_3>t_4>u_5>u_6>0$
\item $u_1>t_1>u_2>t_2>u_3>u_4>t_3>u_5>t_4>u_6>0$
\item $u_1>t_1>u_2>t_2>u_3>u_4>t_3>u_5>u_6>t_4>0$

***

\item $u_1>t_1>u_2>t_2>u_3>u_4>u_5>t_3>t_4>u_6>0$
\item $u_1>t_1>u_2>t_2>u_3>u_4>u_5>t_3>u_6>t_4>0$

***

\item $u_1>t_1>u_2>t_2>u_3>u_4>u_5>u_6>t_3>t_4>0$

***
***

\item $u_1>t_1>u_2>u_3>t_2>t_3>t_4>u_4>u_5>u_6>0$
\item $u_1>t_1>u_2>u_3>t_2>t_3>u_4>t_4>u_5>u_6>0$
\item $u_1>t_1>u_2>u_3>t_2>t_3>u_4>u_5>t_4>u_6>0$
\item $u_1>t_1>u_2>u_3>t_2>t_3>u_4>u_5>u_6>t_4>0$

***

\item $u_1>t_1>u_2>u_3>t_2>u_4>t_3>t_4>u_5>u_6>0$
\item $u_1>t_1>u_2>u_3>t_2>u_4>t_3>u_5>t_4>u_6>0$
\item $u_1>t_1>u_2>u_3>t_2>u_4>t_3>u_5>u_6>t_4>0$

***

\item $u_1>t_1>u_2>u_3>t_2>u_4>u_5>t_3>t_4>u_6>0$
\item $u_1>t_1>u_2>u_3>t_2>u_4>u_5>t_3>u_6>t_4>0$

***

\item $u_1>t_1>u_2>u_3>t_2>u_4>u_5>u_6>t_3>t_4>0$

***
***

\item $u_1>t_1>u_2>u_3>u_4>t_2>t_3>t_4>u_5>u_6>0$
\item $u_1>t_1>u_2>u_3>u_4>t_2>t_3>u_5>t_4>u_6>0$
\item $u_1>t_1>u_2>u_3>u_4>t_2>t_3>u_5>u_6>t_4>0$

***

\item $u_1>t_1>u_2>u_3>u_4>t_2>u_5>t_3>t_4>u_6>0$
\item $u_1>t_1>u_2>u_3>u_4>t_2>u_5>t_3>u_6>t_4>0$

***

\item $u_1>t_1>u_2>u_3>u_4>t_2>u_5>u_6>t_3>t_4>0$

***
***

\item $u_1>t_1>u_2>u_3>u_4>u_5>t_2>t_3>t_4>u_6>0$
\item $u_1>t_1>u_2>u_3>u_4>u_5>t_2>t_3>u_6>t_4>0$

***

\item $u_1>t_1>u_2>u_3>u_4>u_5>t_2>u_6>t_3>t_4>0$

***
***

\item $u_1>t_1>u_2>u_3>u_4>u_5>u_6>t_2>t_3>t_4>0$

***
***
***

\item $u_1>u_2>t_1>t_2>t_3>t_4>u_3>u_4>u_5>u_6>0$
\item $u_1>u_2>t_1>t_2>t_3>u_3>t_4>u_4>u_5>u_6>0$
\item $u_1>u_2>t_1>t_2>t_3>u_3>u_4>t_4>u_5>u_6>0$
\item $u_1>u_2>t_1>t_2>t_3>u_3>u_4>u_5>t_4>u_6>0$
\item $u_1>u_2>t_1>t_2>t_3>u_3>u_4>u_5>u_6>t_4>0$

***

\item $u_1>u_2>t_1>t_2>u_3>t_3>t_4>u_4>u_5>u_6>0$
\item $u_1>u_2>t_1>t_2>u_3>t_3>u_4>t_4>u_5>u_6>0$
\item $u_1>u_2>t_1>t_2>u_3>t_3>u_4>u_5>t_4>u_6>0$
\item $u_1>u_2>t_1>t_2>u_3>t_3>u_4>u_5>u_6>t_4>0$

***

\item $u_1>u_2>t_1>t_2>u_3>u_4>t_3>t_4>u_5>u_6>0$
\item $u_1>u_2>t_1>t_2>u_3>u_4>t_3>u_5>t_4>u_6>0$
\item $u_1>u_2>t_1>t_2>u_3>u_4>t_3>u_5>u_6>t_4>0$

***

\item $u_1>u_2>t_1>t_2>u_3>u_4>u_5>t_3>t_4>u_6>0$
\item $u_1>u_2>t_1>t_2>u_3>u_4>u_5>t_3>u_6>t_4>0$

***

\item $u_1>u_2>t_1>t_2>u_3>u_4>u_5>u_6>t_3>t_4>0$

***
***

\item $u_1>u_2>t_1>u_3>t_2>t_3>t_4>u_4>u_5>u_6>0$
\item $u_1>u_2>t_1>u_3>t_2>t_3>u_4>t_4>u_5>u_6>0$
\item $u_1>u_2>t_1>u_3>t_2>t_3>u_4>u_5>t_4>u_6>0$
\item $u_1>u_2>t_1>u_3>t_2>t_3>u_4>u_5>u_6>t_4>0$

***

\item $u_1>u_2>t_1>u_3>t_2>u_4>t_3>t_4>u_5>u_6>0$
\item $u_1>u_2>t_1>u_3>t_2>u_4>t_3>u_5>t_4>u_6>0$
\item $u_1>u_2>t_1>u_3>t_2>u_4>t_3>u_5>u_6>t_4>0$

***

\item $u_1>u_2>t_1>u_3>t_2>u_4>u_5>t_3>t_4>u_6>0$
\item $u_1>u_2>t_1>u_3>t_2>u_4>u_5>t_3>u_6>t_4>0$

***

\item $u_1>u_2>t_1>u_3>t_2>u_4>u_5>t_3>t_4>u_6>0$

***
***

\item $u_1>u_2>t_1>u_3>u_4>t_2>t_3>t_4>u_5>u_6>0$
\item $u_1>u_2>t_1>u_3>u_4>t_2>t_3>u_5>t_4>u_6>0$
\item $u_1>u_2>t_1>u_3>u_4>t_2>t_3>u_5>u_6>t_4>0$

***

\item $u_1>u_2>t_1>u_3>u_4>t_2>u_5>t_3>t_4>u_6>0$
\item $u_1>u_2>t_1>u_3>u_4>t_2>u_5>t_3>u_6>t_4>0$

***

\item $u_1>u_2>t_1>u_3>u_4>t_2>u_5>u_6>t_3>t_4>0$

***
***

\item $u_1>u_2>t_1>u_3>u_4>u_5>t_2>t_3>t_4>u_6>0$
\item $u_1>u_2>t_1>u_3>u_4>u_5>t_2>t_3>u_6>t_4>0$

***

\item $u_1>u_2>t_1>u_3>u_4>u_5>t_2>u_6>t_3>t_4>0$

***
***

\item $u_1>u_2>t_1>u_3>u_4>u_5>u_6>t_2>t_3>t_4>0$

***
***
***

\item $u_1>u_2>u_3>t_1>t_2>t_3>t_4>u_4>u_5>u_6>0$
\item $u_1>u_2>u_3>t_1>t_2>t_3>u_4>t_4>u_5>u_6>0$
\item $u_1>u_2>u_3>t_1>t_2>t_3>u_4>u_5>t_4>u_6>0$
\item $u_1>u_2>u_3>t_1>t_2>t_3>u_4>u_5>u_6>t_4>0$

***

\item $u_1>u_2>u_3>t_1>t_2>u_4>t_3>t_4>u_5>u_6>0$
\item $u_1>u_2>u_3>t_1>t_2>u_4>t_3>u_5>t_4>u_6>0$
\item $u_1>u_2>u_3>t_1>t_2>u_4>t_3>u_5>u_6>t_4>0$

***

\item $u_1>u_2>u_3>t_1>t_2>u_4>u_5>t_3>t_4>u_6>0$
\item $u_1>u_2>u_3>t_1>t_2>u_4>u_5>t_3>u_6>t_4>0$

***

\item $u_1>u_2>u_3>t_1>t_2>u_4>u_5>u_6>t_3>t_4>0$

***
***

\item $u_1>u_2>u_3>t_1>u_4>t_2>t_3>t_4>u_5>u_6>0$
\item $u_1>u_2>u_3>t_1>u_4>t_2>t_3>u_5>t_4>u_6>0$
\item $u_1>u_2>u_3>t_1>u_4>t_2>t_3>u_5>u_6>t_4>0$

***

\item $u_1>u_2>u_3>t_1>u_4>t_2>u_5>t_3>t_4>u_6>0$
\item $u_1>u_2>u_3>t_1>u_4>t_2>u_5>t_3>u_6>t_4>0$

***

\item $u_1>u_2>u_3>t_1>u_4>t_2>u_5>u_6>t_3>t_4>0$

***
***

\item $u_1>u_2>u_3>t_1>u_4>u_5>t_2>t_3>t_4>u_6>0$
\item $u_1>u_2>u_3>t_1>u_4>u_5>t_2>t_3>u_6>t_4>0$

***

\item $u_1>u_2>u_3>t_1>u_4>u_5>t_2>u_6>t_3>t_4>0$

***
***

\item $u_1>u_2>u_3>t_1>u_4>u_5>u_6>t_2>t_3>t_4>0$

***
***
***

\item $u_1>u_2>u_3>u_4>t_1>t_2>t_3>t_4>u_5>u_6>0$
\item $u_1>u_2>u_3>u_4>t_1>t_2>t_3>u_5>t_4>u_6>0$
\item $u_1>u_2>u_3>u_4>t_1>t_2>t_3>u_5>u_6>t_4>0$

***

\item $u_1>u_2>u_3>u_4>t_1>t_2>u_5>t_3>t_4>u_6>0$
\item $u_1>u_2>u_3>u_4>t_1>t_2>u_5>t_3>u_6>t_4>0$

***

\item $u_1>u_2>u_3>u_4>t_1>t_2>u_5>u_6>t_3>t_4>0$

***
***

\item $u_1>u_2>u_3>u_4>t_1>u_5>t_2>t_3>t_4>u_6>0$
\item $u_1>u_2>u_3>u_4>t_1>u_5>t_2>t_3>u_6>t_4>0$

***

\item $u_1>u_2>u_3>u_4>t_1>u_5>t_2>u_6>t_3>t_4>0$

***
***

\item $u_1>u_2>u_3>u_4>t_1>u_5>u_6>t_2>t_3>t_4>0$

***
***
***

\item $u_1>u_2>u_3>u_4>u_5>t_1>t_2>t_3>t_4>u_6>0$
\item $u_1>u_2>u_3>u_4>u_5>t_1>t_2>t_3>u_6>t_4>0$

***

\item $u_1>u_2>u_3>u_4>u_5>t_1>t_2>u_6>t_3>t_4>0$

***
***

\item $u_1>u_2>u_3>u_4>u_5>t_1>u_6>t_2>t_3>t_4>0$

***
***
***

\item $u_1>u_2>u_3>u_4>u_5>u_6>t_1>t_2>t_3>t_4>0$

\end{enumerate}

\subsection{Shuffle product $\zeta^{(1)}(2,2)\zeta^{(1)}(3,3)$}

Let $DT_{10}$ be the Cartesian product of the standard one dimensional measure on the real line:
\[DT_{10}=dt_1dt_2dt_3dt_4du_1du_2du_3du_4du_5du_6.\]
Let $\Delta$ be the domain of integration obtained by a Cartesian product of a $4$-dimensional simplex with a $6$-dimensional simplex:
\[\Delta=\{(t_1,\dots,t_4,u_1,\dots,u_6)\in \R^{10}\,|\, t_1 > t_2 > t_3 >t_4 > 0  ;  u_1 > u_2 > u_3 > u_4 > u_5 > u_6 > 0\}.\]

\begin{eqnarray*}
&&\zeta^{(1)}(2,2)\zeta^{(1)}(3,3) =\\
& = & \sum_{\alpha_1 \in \mbox {\scriptsize{cone}}} \frac{1}{(\alpha_1)^2 (\alpha_1 + \alpha_2)^2} \sum_{\beta_1 \in \mbox {\scriptsize{cone}}} \frac{1}{(\beta_1)^3 (\beta_1 + \beta_2)^3} \\
& = & \sum \int_\Delta \exp (-\alpha_1t_1 - \alpha_2t_3-\beta_1u_1 - \beta_2u_4)DT_{10} \\
&& =  \zeta^{(1)}(2,2;3,3) + 3\zeta^{(1)}(2,1;4,3) + 3\zeta^{(1)}(2,1;3,4) \\
&& +  3\zeta^{(23)}(2,1;4,3) + 3\zeta^{(23)}(2,1;3,4) \\
&& +  2\zeta^{(23)}(2,2;3,3) + 3\zeta^{(23)}(2,2;2,4) + \zeta^{(23)}(2,3;2,3) \\
&& +  3\zeta^{(23)}(2,3;1,4) + 3\zeta^{(234)}(2,3;1,4) \\
&& +  2\zeta^{(234)}(2,3;2,3) + \zeta^{(234)}(2,3;3,2) + 3\zeta^{(23)}(1,2;4,3) \\
&& +  3\zeta^{(23)}(1,2;3,4) + 4\zeta^{(23)}(1,3;3,3) \\
&& +  6\zeta^{(23)}(1,3;2,4) + 3\zeta^{(23)}(1,4;2,3) + 9\zeta^{(23)}(1,4;1,4) \\
&& +  9\zeta^{(234)}(1,4;1,4) + 6\zeta^{(234)}(1,4;2,3) \\
&& +  3\zeta^{(234)}(1,4;3,2) + 3\zeta^{(234)}(1,3;2,4) + 4\zeta^{(234)}(1,3;3,3) \\
&& +  3\zeta^{(234)}(1,3;4,2) + 3\zeta^{(132)}(1,2;4,3) \\
&& +  3\zeta^{(132)}(1,2;3,4) + 4\zeta^{(132)}(1,3;3,3) + 6\zeta^{(132)}(1,3;2,4) \\
&& +  3\zeta^{(132)}(1,4;2,3) + 9\zeta^{(132)}(1,4;1,4) \\
&& +  9\zeta^{(1342)}(1,4;1,4) + 6\zeta^{(1342)}(1,4;2,3) + 3\zeta^{(1342)}(1,4;3,2) \\
&& +  3\zeta^{(1342)}(1,3;2,4) + 5\zeta^{(1342)}(1,3;3,3) \\
&& +  2\zeta^{(1342)}(1,3;4,2) + 2\zeta^{(132)}(2,2;3,3) + 3\zeta^{(132)}(2,2;2,4) \\
&& +  2\zeta^{(132)}(2,3;2,3) + 6\zeta^{(132)}(2,3;1,4) \\
&& +  6\zeta^{(1342)}(2,3;1,4) + 4\zeta^{(1342)}(2,3;2,3) + 2\zeta^{(1342)}(2,3;3,2) \\
&& +  3\zeta^{(1342)}(2,2;2,4) + 4\zeta^{(1342)}(2,2;3,3) \\
&& +  3\zeta^{(1342)}(3,2;2,3) + 3\zeta^{(132)}(3,2;1,4) + 3\zeta^{(1342)}(3,2;1,4) \\
&& +  \zeta^{(1342)}(3,2;3,2) + 3\zeta^{(1342)}(3,1;2,4) \\
&& +  4\zeta^{(1342)}(3,1;3,3) + 3\zeta^{(1342)}(3,1;4,2) + 3\zeta^{(13)(24)}(3,1;2,4) \\
&& +  4\zeta^{(13)(24)}(3,1;3,3) + 3\zeta^{(13)(24)}(3,1;4,2) \\
&& +  2\zeta^{(13)(24)}(3,2;3,2) + \zeta^{(13)(24)}(3,3;2,2) \\
& = & \zeta^{1}(2,2;3,3) + 3\zeta^{1}(2,1;4,3) + 3\zeta^{1}(2,1;3,4) \\
& + & 3\zeta^{1}(3,1;2,4) + 4\zeta^{1}(3,1;3,3) + 3\zeta^{1}(3,1;4,2) \\
& + & 2\zeta^{1}(3,2;2,3) + 2\zeta^{1}(3,2;3,2) + \zeta^{1}(3,3;2,2) \\
& + & 3\zeta^{(23)}(2,1;4,3) + 3\zeta^{(23)}(2,1;3,4) + 8\zeta^{(23)}(2,2;3,3) \\
& + & 9\zeta^{(23)}(2,2;2,4) + 9\zeta^{(23)}(2,3;2,3) + 18\zeta^{(23)}(2,3;1,4) \\
& + & 3\zeta^{(23)}(2,3;3,2) + 3\zeta^{(23)}(2,2;4,2) + 6\zeta^{(23)}(1,2;4,3) \\
& + & 6\zeta^{(23)}(1,2;3,4) + 8\zeta^{(23)}(1,3;3,3) + 18\zeta^{(23)}(1,3;2,4) \\
& + & 18\zeta^{(23)}(1,4;2,3) + 36\zeta^{(23)}(1,4;1,4) + 6\zeta^{(23)}(1,4;3,2) \\
& + & 9\zeta^{(23)}(1,3;3,3) + 5\zeta^{(23)}(1,3;4,2) + 3\zeta^{(23)}(3,2;2,3) \\
& + & 6\zeta^{(23)}(3,2;1,4) + \zeta^{(23)}(3,2;3,2) + 3\zeta^{(23)}(3,1;2,4) \\
& + & 4\zeta^{(23)}(3,1;3,3) + 3\zeta^{(23)}(3,1;4,2) . 
\end{eqnarray*}

\subsection{Shuffle  product $\zeta^{(1)}(2,2)\zeta^{(1)}(2,4)$}

\begin{eqnarray*}
&&\zeta^{(1)}(2,2)\zeta^{(1)}(2,4) =\\
& = & \sum_{\alpha_1 \in \mbox {\scriptsize{cone}}} \frac{1}{(\alpha_1)^2 (\alpha_1 + \alpha_2)^2} \sum_{\beta_1 \in \mbox {\scriptsize{cone}}} \frac{1}{(\beta_1)^2 (\beta_1 + \beta_2)^4} \\
& = & \sum \int_\Delta \exp (-\alpha_1t_1 - \alpha_2t_3-\beta_1u_1 - \beta_2u_3)DT_{10} \\
& =  &\zeta^{(1)}(2,2;2,4) + 2\zeta^{(1)}(2,1;3,4) + 4\zeta^{(1)}(2,1;2,5) \\
&& +  2\zeta^{(23)}(2,1;3,4) + 4\zeta^{(23)}(2,1;2,5) + \zeta^{(23)}(2,2;2,4) \\
&& +  4\zeta^{(23)}(2,2;1,5) + 4\zeta^{(234)}(2,2;1,5) + 3\zeta^{(234)}(2,2;2,4) \\
&& +  2\zeta^{(234)}(2,2;3,3) + \zeta^{(234)}(2,2;4,2) + 2\zeta^{(23)}(1,2;3,4) \\
&& +  4\zeta^{(23)}(1,2;2,5) + 2\zeta^{(23)}(1,3;2,4) + 8\zeta^{(23)}(1,3;1,5) \\
&& +  8\zeta^{(234)}(1,3;1,5) + 6\zeta^{(234)}(1,3;2,4) + 4\zeta^{(234)}(1,3;3,3) \\
&& +  2\zeta^{(234)}(1,3;4,2) + 4\zeta^{(234)}(1,2;2,5) + 6\zeta^{(234)}(1,2;3,4) \\
&& +  6\zeta^{(234)}(1,2;4,3) + 4\zeta^{(234)}(1,2;5,2) + 2\zeta^{(132)}(1,2;3,4) \\
&& +  4\zeta^{(132)}(1,2;2,5) + 2\zeta^{(132)}(1,3;2,4) + 8\zeta^{(132)}(1,3;1,5) \\
&& +  8\zeta^{(1342)}(1,3;1,5) + 6\zeta^{(1342)}(1,3;2,4) + 4\zeta^{(1342)}(1,3;3,3) \\
&& +  2\zeta^{(1342)}(1,3;4,2) + 4\zeta^{(1342)}(1,2;2,5) + 6\zeta^{(1342)}(1,2;3,4) \\
&& +  6\zeta^{(1342)}(1,2;4,3) + 4\zeta^{(1342)}(1,2;5,2) + \zeta^{(132)}(2,2;2,4) \\
&& +  4\zeta^{(132)}(2,2;1,5) + 4\zeta^{(1342)}(2,2;1,5) + 3\zeta^{(1342)}(2,2;2,4) \\
&& +  2\zeta^{(1342)}(2,2;3,3) + \zeta^{(1342)}(2,2;4,2) + 4\zeta^{(1342)}(2,1;2,5) \\
&& +  6\zeta^{(1342)}(2,1;3,4) + 6\zeta^{(1342)}(2,1;4,3) + 4\zeta^{(1342)}(2,1;5,2) \\
&& +  4\zeta^{(13)(24)}(2,1;2,5) + 6\zeta^{(13)(24)}(2,1;3,4) + 6\zeta^{(13)(24)}(2,1;4,3) \\
&& +  4\zeta^{(13)(24)}(2,1;5,2) + 3\zeta^{(13)(24)}(2,2;2,4) + 4\zeta^{(13)(24)}(2,2;3,3) \\
&& +  3\zeta^{(13)(24)}(2,2;4,2) + 2\zeta^{(13)(24)}(2,3;2,3) + 2\zeta^{(13)(24)}(2,3;3,2) \\
&& +  \zeta^{(13)(24)}(2,4;2,2) \\
& = & 4\zeta^{1}(2,2;2,4) + 8\zeta^{1}(2,1;3,4) + 8\zeta^{1}(2,1;2,5) \\
& + & 6\zeta^{1}(2,1;4,3) + 4\zeta^{1}(2,1;5,2) + 4\zeta^{1}(2,2;3,3) \\
& + & 3\zeta^{1}(2,2;4,2) + 2\zeta^{1}(2,3;2,3) + 2\zeta^{1}(2,3;3,2) \\
& + & \zeta^{1}(2,4;2,2) + 8\zeta^{(23)}(2,1;3,4) + 8\zeta^{(23)}(2,1;2,5) \\
& + & 8\zeta^{(23)}(2,2;2,4) + 16\zeta^{(23)}(2,2;1,5) + 4\zeta^{(23)}(2,2;3,3) \\
& + & 2\zeta^{(23)}(2,2;4,2) + 16\zeta^{(23)}(1,2;3,4) + 16\zeta^{(23)}(1,2;2,5) \\
& + & 16\zeta^{(23)}(1,3;2,4) + 32\zeta^{(23)}(1,3;1,5) \\
& + & 8\zeta^{(23)}(1,3;3,3) + 4\zeta^{(23)}(1,3;4,2) + 12\zeta^{(23)}(1,2;4,3) \\
& + & 8\zeta^{(23)}(1,2;5,2) + 6\zeta^{(23)}(2,1;4,3) + 4\zeta^{(23)}(2,1;5,2) .
\end{eqnarray*}

\subsection{Shuffle  product $\zeta^{(1)}(2,2)\zeta^{(1)}(1,5)$}

\begin{eqnarray*}
&&\zeta^{(1)}(2,2)\zeta^{(1)}(1,5) =\\
& = & \sum_{\alpha_1 \in \mbox {\scriptsize{cone}}} \frac{1}{(\alpha_1)^2 (\alpha_1 + \alpha_2)^2} \sum_{\beta_1 \in \mbox {\scriptsize{cone}}} \frac{1}{\beta_1 (\beta_1 + \beta_2)^4} \\
& = & \sum \int_\Delta \exp (-\alpha_1t_1 - \alpha_2t_3-\beta_1u_1 - \beta_2u_2)DT_{10} \\
& = & \zeta^{(1)}(2,2;1,5) + \zeta^{(1)}(2,1;2,5) + 5\zeta^{(1)}(2,1;1,6) \\
&& +  \zeta^{(23)}(2,1;2,5) + 5\zeta^{(23)}(2,1;1,6) + 5\zeta^{(234)}(2,1;1,6) \\
&& +  4\zeta^{(234)}(2,1;2,5) + 3\zeta^{(234)}(2,1;3,4) + 2\zeta^{(234)}(2,1;4,3) \\
&& +  \zeta^{(234)}(2,1;5,2) + \zeta^{(23)}(1,2;2,5) + 5\zeta^{(23)}(1,2;1,6) \\ 
&& +  5\zeta^{(234)}(1,2;1,6) + 4\zeta^{(234)}(1,2;2,5) + 3\zeta^{(234)}(1,2;3,4) \\
&& +  2\zeta^{(234)}(1,2;4,3) + \zeta^{(234)}(1,2;5,2) + 5\zeta^{(234)}(1,1;2,6) \\
&& +  8\zeta^{(234)}(1,1;3,5) + 9\zeta^{(234)}(1,1;4,4) + 8\zeta^{(234)}(1,1;5,3) \\
&& +  5\zeta^{(234)}(1,1;6,2) + \zeta^{(132)}(1,2;2,5) + 5\zeta^{(132)}(1,2;1,6) \\
&& +  5\zeta^{(1342)}(1,2;1,6) + 4\zeta^{(1342)}(1,2;2,5) + 3\zeta^{(1342)}(1,2;3,4) \\
&& +  2\zeta^{(1342)}(1,2;4,3) + \zeta^{(1342)}(1,2;5,2) + 5\zeta^{(1342)}(1,1;2,6) \\
&& +  8\zeta^{(1342)}(1,1;3,5) + 9\zeta^{(1342)}(1,1;4,4) + 8\zeta^{(1342)}(1,1;5,3) \\
&& +  5\zeta^{(1342)}(1,1;6,2) + 5\zeta^{(13)(24)}(1,1;2,6) + 8\zeta^{(13)(24)}(1,1;3,5) \\
&& +  9\zeta^{(13)(24)}(1,1;4,4) + 8\zeta^{(13)(24)}(1,1;5,3) + 5\zeta^{(13)(24)}(1,1;6,2) \\
&& +  4\zeta^{(13)(24)}(1,2;2,5) + 6\zeta^{(13)(24)}(1,2;3,4) + 6\zeta^{(13)(24)}(1,2;4,3) \\
&& +  4\zeta^{(13)(24)}(1,2;5,2) + 3\zeta^{(13)(24)}(1,3;2,4) + 4\zeta^{(13)(24)}(1,3;3,3) \\
&& +  3\zeta^{(13)(24)}(1,3;4,2) + 2\zeta^{(13)(24)}(1,4;2,3) + 2\zeta^{(13)(24)}(1,4;3,2) \\
&& +  \zeta^{(13)(24)}(1,5;2,2) \\
& = & \zeta^{1}(2,2;1,5) + \zeta^{1}(2,1;2,5) + 5\zeta^{1}(2,1;1,6) \\
& + & 5\zeta^{1}(1,1;2,6) + 8\zeta^{1}(1,1;3,5) + 9\zeta^{1}(1,1;4,4) \\
& + & 8\zeta^{1}(1,1;5,3) + 5\zeta^{1}(1,1;6,2) + 4\zeta^{1}(1,2;2,5) \\
& + & 6\zeta^{1}(1,2;3,4) + 6\zeta^{1}(1,2;4,3) + 4\zeta^{1}(1,2;5,2) \\
& + & 3\zeta^{1}(1,3;2,4) + 4\zeta^{1}(1,3;3,3) + 3\zeta^{1}(1,3;4,2) \\
& + & 2\zeta^{1}(1,4;2,3) + 2\zeta^{1}(1,4;3,2) + \zeta^{1}(1,5;2,2) \\
& + & 5\zeta^{(23)}(2,1;2,5) + 10\zeta^{(23)}(2,1;1,6) + 3\zeta^{(23)}(1,3;4,2) \\
& + & 2\zeta^{(23)}(2,1;4,3) + \zeta^{(23)}(2,1;5,2) + 10\zeta^{(23)}(1,2;2,5) \\
& + & 20\zeta^{(23)}(1,2;1,6) + 6\zeta^{(23)}(1,2;3,4) + 4\zeta^{(23)}(1,2;4,3) \\
& + & 2\zeta^{(23)}(1,2;5,2) + 10\zeta^{(23)}(1,1;2,6) + 16\zeta^{(23)}(1,1;3,5) \\
& + & 18\zeta^{(23)}(1,1;4,4) + 16\zeta^{(23)}(1,1;5,3) + 10\zeta^{(23)}(1,1;6,2) .
\end{eqnarray*}

\subsection{Shuffle  product  $\zeta^{(1)}(1,3)\zeta^{(1)}(3,3)$}

\begin{eqnarray*}
&&\zeta^{(1)}(1,3)\zeta^{(1)}(3,3) =\\
& = & \sum_{\alpha_1 \in \mbox {\scriptsize{cone}}} \frac{1}{\alpha_1 (\alpha_1 + \alpha_2)^3} \sum_{\beta_1 \in \mbox {\scriptsize{cone}}} \frac{1}{(\beta_1)^3 (\beta_1 + \beta_2)^3} \\
& = & \sum \int_\Delta \exp (-\alpha_1t_1 - \alpha_2t_2-\beta_1u_1 - \beta_2u_4)DT_{10} \\
& = & \zeta^{(1)}(1,3;3,3) + 3\zeta^{(1)}(1,2;4,3) + 3\zeta^{(1)}(1,2;3,4) \\ 
&& +  6\zeta^{(1)}(1,1;5,3) + 9\zeta^{(1)}(1,1;4,4) + 6\zeta^{(1)}(1,1;3,5) \\
&& +  6\zeta^{(23)}(1,1;5,3) + 9\zeta^{(23)}(1,1;4,4) + 6\zeta^{(23)}(1,1;3,5) \\
&& +  3\zeta^{(23)}(1,2;4,3) + 6\zeta^{(23)}(1,2;3,4) + 6\zeta^{(23)}(1,2;2,5) \\
&& +  \zeta^{(23)}(1,3;3,3) + 3\zeta^{(23)}(1,3;2,4) + 6\zeta^{(23)}(1,3;1,5) \\
&& +  6\zeta^{(234)}(1,3;1,5) + \zeta^{(234)}(1,3;3,3) + 6\zeta^{(132)}(1,1;5,3) \\
&& +  9\zeta^{(132)}(1,1;4,4) + 6\zeta^{(132)}(1,1;3,5) + 3\zeta^{(132)}(1,2;4,3) \\
&& +  6\zeta^{(132)}(1,2;3,4) + 6\zeta^{(132)}(1,2;2,5) + \zeta^{(132)}(1,3;3,3) \\
&& +  3\zeta^{(132)}(1,3;2,4) + 6\zeta^{(132)}(1,3;1,5) + 6\zeta^{(1342)}(1,3;1,5) \\
&& +  3\zeta^{(1342)}(1,3;2,4) + \zeta^{(1342)}(1,3;3,3) + 3\zeta^{(132)}(2,1;4,3) \\
&& +  6\zeta^{(132)}(2,1;3,4) + 6\zeta^{(132)}(2,1;2,5) + \zeta^{(132)}(2,2;3,3) \\
&& +  3\zeta^{(132)}(2,2;2,4) + 6\zeta^{(132)}(2,2;1,5) + 6\zeta^{(1342)}(2,2;1,5) \\
&& +  3\zeta^{(1342)}(2,2;2,4) + \zeta^{(1342)}(2,2;3,3) + \zeta^{(1342)}(3,1;3,3) \\
&& +  3\zeta^{(132)}(3,1;2,4) + 6\zeta^{(132)}(3,1;1,5) + 6\zeta^{(1342)}(3,1;1,5) \\
&& +  3\zeta^{(1342)}(3,1;2,4) + \zeta^{(1342)}(3,1;3,3) + 6\zeta^{(13)(24)}(3,1;1,5) \\
&& +  3\zeta^{(13)(24)}(3,1;2,4) + \zeta^{(13)(24)}(3,1;3,3) + 3\zeta^{(13)(24)}(3,2;1,4) \\
&& +  \zeta^{(13)(24)}(3,2;2,3) + \zeta^{(13)(24)}(3,3;1,3) \\
& = & \zeta^{1}(1,3;3,3) + 3\zeta^{1}(1,2;4,3) + 3\zeta^{1}(1,2;3,4) \\
& + & 6\zeta^{1}(1,1;5,3) + 9\zeta^{1}(1,1;4,4) + 6\zeta^{1}(1,1;3,5) \\
& + & 6\zeta^{1}(3,1;1,5) + 3\zeta^{1}(3,1;2,4) + \zeta^{1}(3,1;3,3) \\
& + & 3\zeta^{1}(3,2;1,4) + \zeta^{1}(3,2;2,3) + \zeta^{1}(3,3;1,3) \\
& + & 12\zeta^{(23)}(1,1;5,3) + 18\zeta^{(23)}(1,1;4,4) + 12\zeta^{(23)}(1,1;3,5) \\
& + & 6\zeta^{(23)}(1,2;4,3) + 12\zeta^{(23)}(1,2;3,4) + 12\zeta^{(23)}(1,2;2,5) \\
& + & 4\zeta^{(23)}(1,3;3,3) + 12\zeta^{(23)}(1,3;2,4) + 24\zeta^{(23)}(1,3;1,5) \\
& + & 3\zeta^{(23)}(2,1;4,3) + 6\zeta^{(23)}(2,1;3,4) + 6\zeta^{(23)}(2,1;2,5) + 2\zeta^{(23)}(2,2;3,3) \\
& + & 6\zeta^{(23)}(2,2;2,4) + 12\zeta^{(23)}(2,2;1,5) + 2\zeta^{(23)}(3,1;3,3) \\
& + & 6\zeta^{(23)}(3,1;2,4) + 12\zeta^{(23)}(3,1;1,5) .
\end{eqnarray*}

\subsection{Shuffle  product $\zeta^{(1)}(1,3)\zeta^{(1)}(2,4)$}

\begin{eqnarray*}
&&\zeta^{(1)}(1,3)\zeta^{(1)}(2,4) =\\
& = & \sum_{\alpha_1 \in \mbox {\scriptsize{cone}}} \frac{1}{\alpha_1 (\alpha_1 + \alpha_2)^3} \sum_{\beta_1 \in \mbox {\scriptsize{cone}}} \frac{1}{(\beta_1)^2 (\beta_1 + \beta_2)^4} \\
& = & \sum \int_\Delta \exp (-\alpha_1t_1 - \alpha_2t_2-\beta_1u_1 - \beta_2u_3)DT_{10} \\
& = & \zeta^{(1)}(1,3;2,4) + 2\zeta^{(1)}(1,2;3,4) + 4\zeta^{(1)}(1,2;2,5) \\
&& +  3\zeta^{(1)}(1,1;4,4) + 8\zeta^{(1)}(1,1;3,5) + 10\zeta^{(1)}(1,1;2,6) \\ 
&& +  3\zeta^{(23)}(1,1;4,4) + 8\zeta^{(23)}(1,1;3,5) + 10\zeta^{(23)}(1,1;2,6) \\
&& +  \zeta^{(23)}(1,2;3,4) + 4\zeta^{(23)}(1,2;2,5) + 10\zeta^{(23)}(1,2;1,6) \\
&& +  10\zeta^{(234)}(1,2;1,6) + 6\zeta^{(234)}(1,2;2,5) + 3\zeta^{(234)}(1,2;3,4) \\
&& +  \zeta^{(234)}(1,2;4,3) + 3\zeta^{(132)}(1,1;4,4) + 8\zeta^{(132)}(1,1;3,5) \\
&& +  10\zeta^{(132)}(1,1;2,6) + \zeta^{(132)}(1,2;3,4) + 4\zeta^{(132)}(1,2;2,5) \\
&& +  10\zeta^{(132)}(1,2;1,6) + 10\zeta^{(1342)}(1,2;1,6) + 6\zeta^{(1342)}(1,2;2,5) \\
&& +  3\zeta^{(1342)}(1,2;3,4) + \zeta^{(1342)}(1,2;4,3) + \zeta^{(132)}(2,1;3,4) \\
&& +  4\zeta^{(132)}(2,1;2,5) + 10\zeta^{(132)}(2,1;1,6) + 10\zeta^{(1342)}(2,1;1,6) \\ 
&& +  6\zeta^{(1342)}(2,1;2,5) + 3\zeta^{(1342)}(2,1;3,4) + \zeta^{(1342)}(2,1;4,3) \\
&& +  10\zeta^{(13)(24)}(2,1;1,6) + 6\zeta^{(13)(24)}(2,1;2,5) + 3\zeta^{(13)(24)}(2,1;3,4) \\
&& +  \zeta^{(13)(24)}(2,1;4,3) + 6\zeta^{(13)(24)}(2,2;1,5) + 3\zeta^{(13)(24)}(2,2;2,4) \\
&& +  \zeta^{(13)(24)}(2,2;3,3) + 3\zeta^{(13)(24)}(2,3;1,4) + \zeta^{(13)(24)}(2,3;2,3) \\
&& +  \zeta^{(13)(24)}(2,4;1,3) \\
& = & \zeta^{1}(1,3;2,4) + 2\zeta^{1}(1,2;3,4) + 4\zeta^{1}(1,2;2,5) \\
& + & 3\zeta^{1}(1,1;4,4) + 8\zeta^{1}(1,1;3,5) + 10\zeta^{1}(1,1;2,6) \\
& + & 10\zeta^{1}(2,1;1,6) + 6\zeta^{1}(2,1;2,5) + 3\zeta^{1}(2,1;3,4) \\
& + & \zeta^{1}(2,1;4,3) + 6\zeta^{1}(2,2;1,5) + 3\zeta^{1}(2,2;2,4) \\
& + & \zeta^{1}(2,2;3,3) + 3\zeta^{1}(2,3;1,4) + \zeta^{1}(2,3;2,3) \\
& + & \zeta^{1}(2,4;1,3) + 6\zeta^{(23)}(1,1;4,4) + 16\zeta^{(23)}(1,1;3,5) \\
& + & 20\zeta^{(23)}(1,1;2,6) + 8\zeta^{(23)}(1,2;3,4) + 20\zeta^{(23)}(1,2;2,5) \\
& + & 2\zeta^{(23)}(1,2;4,3) + 4\zeta^{(23)}(2,1;3,4) + 10\zeta^{(23)}(2,1;2,5) \\
& + & 20\zeta^{(23)}(2,1;1,6) + \zeta^{(23)}(2,1;4,3) + 40\zeta^{(23)}(1,2;1,6).
\end{eqnarray*}

\subsection{Shuffle  product $\zeta^{(1)}(1,3)\zeta^{(1)}(1,5)$}

\begin{eqnarray*}
&&\zeta^{(1)}(1,3)\zeta^{(1)}(1,5) =\\
& = & \sum_{\alpha_1 \in \mbox {\scriptsize{cone}}} \frac{1}{\alpha_1 (\alpha_1 + \alpha_2)^3} \sum_{\beta_1 \in \mbox {\scriptsize{cone}}} \frac{1}{\beta_1 (\beta_1 + \beta_2)^5} \\
& = & \sum \int_\Delta \exp (-\alpha_1t_1 - \alpha_2t_2-\beta_1u_1 - \beta_2u_2)DT_{10} \\
& = & \zeta^{(1)}(1,3;1,5) + \zeta^{(1)}(1,2;2,5) + 5\zeta^{(1)}(1,2;1,6) \\
&& +  \zeta^{(1)}(1,1;3,5) + 5\zeta^{(1)}(1,1;2,6) + 15\zeta^{(1)}(1,1;1,7) \\
&& +  \zeta^{(23)}(1,1;3,5) + 5\zeta^{23)}(1,1;2,6) + 15\zeta^{(23)}(1,1;1,7) \\
&& +  15\zeta^{(234)}(1,1;1,7) + 10\zeta^{(234)}(1,1;2,6) + 6\zeta^{(234)}(1,1;3,5) \\
&& +  3\zeta^{(234)}(1,1;4,4) + \zeta^{(234)}(1,1;5,3) + \zeta^{(132)}(1,1;3,5) \\
&& +  5\zeta^{(132)}(1,1;2,6) + 15\zeta^{(132)}(1,1;1,7) + 15\zeta^{(1342)}(1,1;1,7) \\
&& +  10\zeta^{(1342)}(1,1;2,6) + 6\zeta^{(1342)}(1,1;3,5) + 3\zeta^{(1342)}(1,1;4,4) \\
&& +  \zeta^{(1342)}(1,1;5,3) + 15\zeta^{(13)(24)}(1,1;1,7) + 10\zeta^{(13)(24)}(1,1;2,6) \\
&& +  6\zeta^{(13)(24)}(1,1;3,5) + 3\zeta^{(13)(24)}(1,1;4,4) + \zeta^{(13)(24)}(1,1;5,3) \\
&& +  10\zeta^{(13)(24)}(1,2;1,6) + 6\zeta^{(13)(24)}(1,2;2,5) + 3\zeta^{(13)(24)}(1,2;3,4) \\
&& +  \zeta^{(13)(24)}(1,2;4,3) + 6\zeta^{(13)(24)}(1,3;1,5) + 3\zeta^{(13)(24)}(1,3;2,4) \\
&& +  \zeta^{(13)(24)}(1,3;3,3) + 3\zeta^{(13)(24)}(1,4;1,4) + \zeta^{(13)(24)}(1,4;2,3) + \zeta^{(13)(24)}(1,5;1,3) \\
& = & 7\zeta^{1}(1,3;1,5) + 7\zeta^{1}(1,2;2,5) + 15\zeta^{1}(1,2;1,6) \\
& + & 7\zeta^{1}(1,1;3,5) + 15\zeta^{1}(1,1;2,6) + 30\zeta^{1}(1,1;1,7) \\
& + & 3\zeta^{1}(1,1;4,4) + \zeta^{1}(1,1;5,3) + 3\zeta^{1}(1,2;3,4) \\
& + & \zeta^{1}(1,2;4,3) + 3\zeta^{1}(1,3;2,4) + \zeta^{1}(1,3;3,3) \\
& + & 3\zeta^{1}(1,4;14) + \zeta^{1}(1,4;2,3) + \zeta^{1}(1,5;1,3) \\
& + & 14\zeta^{(23)}(1,1;3,5) + 30\zeta^{(23)}(1,1;2,6) + 60\zeta^{(23)}(1,1;1,7) \\
& + & 6\zeta^{(23)}(1,1;4,4) + 2\zeta^{(23)}(1,1;5,3).
\end{eqnarray*}

\subsection{Formula for the shuffle product of $\zeta_{K,C}(2)\zeta_{K,C}(3)$}
Using the self shuffles of  $\zeta_{K,C}(2)$ and of $\zeta_{K,C}(3)$ together with the simplified formulas for the previous six subsections, assuming Galois invariance of the cone, we obtain the following Theorem. 

\begin{theorem}
The shuffle product of $\zeta_{K,C}(2)$ times $\zeta_{K,C}(3)$ in terms of refinements of MDZVs is given by
\begin{eqnarray*}
&&\zeta_{K,C}(2)\zeta_{K,C}(3) =\\
& = & 152\zeta^{1}(1,3;3,3) + 216\zeta^{1}(1,2;4,3) + 360\zeta^{1}(1,2;3,4) \\
&& +  288\zeta^{1}(1,1;5,3) + 504\zeta^{1}(1,1;4,4) + 768\zeta^{1}(1,1;3,5) \\
&& +  48\zeta^{1}(3,1;1,5) + 36\zeta^{1}(3,1;2,4) + 24\zeta^{1}(3,1;3,3) \\
&& +  24\zeta^{1}(3,2;1,4) + 16\zeta^{1}(3,2;2,3) + 8\zeta^{1}(3,3;1,3) \\
&& +  576\zeta^{(23)}(1,1;5,3) + 1008\zeta^{(23)}(1,1;4,4) + 1536\zeta^{(23)}(1,1;3,5) \\ 
&& +  360\zeta^{(23)}(1,2;4,3) + 648\zeta^{(23)}(1,2;3,4) + 1008\zeta^{(23)}(1,2;2,5) \\
&& +  196\zeta^{(23)}(1,3;3,3) + 360\zeta^{(23)}(1,3;2,4) + 576\zeta^{(23)}(1,3;1,5) \\
&& +  180\zeta^{(23)}(2,1;4,3) + 252\zeta^{(23)}(2,1;3,4) + 504\zeta^{(23)}(2,1;2,5) \\
&& +  96\zeta^{(23)}(2,2;3,3) + 180\zeta^{(23)}(2,2;2,4) + 288\zeta^{(23)}(2,2;1,5) \\
&& +  32\zeta^{(23)}(3,1;3,3) + 60\zeta^{(23)}(3,1;2,4) + 96\zeta^{(23)}(3,1;1,5) \\
&& +  240\zeta^{1}(1,3;2,4) + 528\zeta^{1}(1,2;2,5) + 1080\zeta^{1}(1,1;2,6) \\
&& +  360\zeta^{1}(2,1;1,6) + 264\zeta^{1}(2,1;2,5) + 180\zeta^{1}(2,1;3,4) \\
&& +  108\zeta^{1}(2,1;4,3) + 168\zeta^{1}(2,2;1,5) + 120\zeta^{1}(2,2;2,4) + 76\zeta^{1}(2,2;3,3) \\
&& +  72\zeta^{1}(2,3;1,4) + 48\zeta^{1}(2,3;2,3) + 24\zeta^{1}(2,4;1,3) \\
&& +  2160\zeta^{(23)}(1,1;2,6) + 720\zeta^{(23)}(2,1;1,6) + 1440\zeta^{(23)}(1,2;1,6) \\
&& +  336\zeta^{1}(1,3;1,5) + 720\zeta^{1}(1,2;1,6) + 1440\zeta^{1}(1,1;1,7) \\
&& +  144\zeta^{1}(1,4;1,4) + 96\zeta^{1}(1,4;2,3) + 48\zeta^{1}(1,5;1,3) \\
&& +  2880\zeta^{(23)}(1,1;1,7) + 12\zeta^{1}(3,1;4,2) + 8\zeta^{1}(3,2;3,2) \\
&& + 4\zeta^{1}(3,3;2,2) + 36\zeta^{(23)}(2,3;2,3) + 72\zeta^{(23)}(2,3;1,4) \\
&& + 12\zeta^{(23)}(2,3;3,2) + 36\zeta^{(23)}(2,2;4,2) + 72\zeta^{(23)}(1,4;2,3) \\
&& +  144\zeta^{(23)}(1,4;1,4) + 24\zeta^{(23)}(1,4;3,2) + 140\zeta^{(23)}(1,3;4,2) \\
&& +  12\zeta^{(23)}(3,2;2,3) + 24\zeta^{(23)}(3,2;1,4) + 4\zeta^{(23)}(3,2;3,2) \\
&& +  12\zeta^{(23)}(3,1;4,2) + 48\zeta^{1}(2,1;5,2) + 36\zeta^{1}(2,2;4,2) \\
&& +  24\zeta^{1}(2,3;3,2) + 12\zeta^{1}(2,4;2,2) + 144\zeta^{(23)}(1,2;5,2) \\
&& +  72\zeta^{(23)}(2,1;5,2) + 120\zeta^{1}(1,1;6,2) + 96\zeta^{1}(1,2;5,2) \\
&& +  72\zeta^{1}(1,3;4,2) + 48\zeta^{1}(1,4;3,2) + 24\zeta^{1}(1,5;2,2) + 240\zeta^{(23)}(1,1;6,2) . 
\end{eqnarray*}
\end{theorem}

\renewcommand{\em}{\textrm}

\begin{small}

\renewcommand{\refname}{ {\flushleft\normalsize\bf{References}} }
    
\end{small}

\vspace{.4cm}
Michael Dotzel,
University of Missouri Columbia,
230 Jesse Hall,
Columbia, MO 65211,
USA


\vspace{.2cm}
Ivan Horozov,
City University of New York, Bronx,
Department of Mathematics and Computer Science,
CP 315,
2155 University Ave,
Bronx, NY 10453,
USA

\end{document}